\documentclass[12pt]{article}

\usepackage{amsmath,amssymb,amsthm}

\usepackage{upref}
\usepackage{latexsym}
\newcommand{\beqn}{\begin{eqnarray}}
 \newcommand{\eeqn}{\end{eqnarray}}

 \newcommand{\be}{\begin{equation}}
 \newcommand{\ee}{\end{equation}}
 \newcommand{\ba}{\begin{array}}
 \newcommand{\ea}{\end{array}}
\newcommand{\pa}{\partial}
 \newcommand{\re}{\ref}
 \newcommand{\ci}{\cite}
 \newcommand{\ds}{\displaystyle}
 \newcommand{\la}{\label}
\newcommand{\rIm}{{\rm Im\5}}
 \newcommand{\rRe}{{\rm Re\5}}
 \newcommand{\fr}{\frac}

\newcommand{\ov}{\overline}

\newcommand{\loota}{\hbox{\enspace{\vrule height 7pt depth 0pt width
      7pt}}}
\newcommand{\bo}{{\hfill\loota}}

\renewcommand{\Pr}{{\bf Proof~}}

\newcommand{\bB}{{\bf B}}
\newcommand{\bC}{{\bf C}}
\newcommand{\bP}{{\bf P}}
\newcommand{\bR}{{\bf R}}
\newcommand{\bF}{{\bf F}}
\newcommand{\cO}{{\cal O}}
\newcommand{\cS}{{\cal S}}
\newcommand{\cC}{{\cal C}}
\newcommand{\ve}{\varepsilon}

\newcommand{\de}{\delta}
\newcommand{\al}{\alpha}

\newcommand{\om}{\omega}
\newcommand{\Om}{\Omega}
\newcommand{\na}{\nabla}
\newcommand{\lam}{\lambda}

\newcommand{\5}{{\hspace{0.5mm}}}
\newcommand{\R}{\mathbb{R}}
\newcommand{\C}{\mathbb{C}}
\newcommand{\Z}{\mathbb{Z}}

\renewcommand{\theequation}{\thesection.\arabic{equation}}
\renewcommand{\thesection}{\arabic{section}}
\renewcommand{\thesubsection}{\arabic{section}.\arabic{subsection}}
\newtheorem{theorem}{Theorem}[section]
\renewcommand{\thetheorem}{\arabic{section}.\arabic{theorem}}
\newtheorem{defin}[theorem]{Definition}
\newtheorem{lemma}[theorem]{Lemma}

\newtheorem{remark}[theorem]{Remark}
\newtheorem{remarks}[theorem]{Remarks}
\newtheorem{cor}[theorem]{Corollary}
\newtheorem{pro}[theorem]{Proposition}

\newcommand{\bd}{\begin{defin}}
 \newcommand{\ed}{\end{defin}}
\newcommand{\bt}{\begin{theorem}}
 \newcommand{\et}{\end{theorem}}
\newcommand{\bp}{\begin{pro}}
 \newcommand{\ep}{\end{pro}}
\newcommand{\bl}{\begin{lemma}}
 \newcommand{\el}{\end{lemma}}
\newcommand{\bc}{\begin{cor}}
 \newcommand{\ec}{\end{cor}}
\newcommand{\br}{\begin{remark} }
 \newcommand{\er}{\end{remark}}
\newcommand{\brs}{\begin{remarks}}
 \newcommand{\ers}{\end{remarks}}

\begin{document}

\begin{titlepage}

\bigskip\bigskip\bigskip

\begin{center}
 {\Large\bf
On Asymptotic Stability of Solitary Waves
\bigskip\\
in  Discrete Schr\"odinger Equation
\bigskip\\
Coupled to Nonlinear Oscillator}
\vspace{1cm}
\\
{\large E.~A.~Kopylova}
\footnote{Supported partly by  FWF grant P19138-N13,
DFG grant 436 RUS 113/929/0-1
and RFBR grant 07-01-00018a.}\\
{\it Institute for Information Transmission Problems RAS\\
B.Karetnyi 19, Moscow 101447,GSP-4, Russia}\\
e-mail:~ek@vpti.vladimir.ru
\end{center}

\date{}

\markboth{E.Kopylova}
  {On Asymptotic Stability of Solitary Waves}

%%%%%%%%%%%%%%%%%%%%%%%%%%%%%%%%%%%%%%%%%%%%%%%%%%%%%%%%%%%%%%%%%%%%%%%%%%%%%%%%%%
\vspace{0.5cm}
\begin{abstract}
The long-time asymptotics is analyzed for finite energy solutions of the
1D  discrete Schr\"odinger equation  coupled to a nonlinear oscillator.
The coupled system is invariant with respect to the phase
rotation group $U(1)$.
For initial states close to a solitary wave, the solution converges to a sum
of another solitary wave and dispersive wave which is a solution to the free
Schr\"odinger equation.
The proofs use the strategy of Buslaev-Perelman \ci{BP}:
the linerization of the dynamics on the solitary manifold,
the symplectic orthogonal projection, method of majorants, etc.
\end{abstract}

\end{titlepage}
%%%%%%%%%%%%%%%%%%%%%%%%%%%%%%%%%%%%%%%%%%%%%%%%%%%%%%%%%%%%%%%%%%%%%%%%%%%%%%%%%%
%%%%%%%%%%%%%%%%%%%%%%%%%%%%%%%%%%%%%%%%%%%%%%%%%%%%%%%%%%%%%%%%%%%%%%%%%%%%%%%%%
\setcounter{equation}{0}
\section{Introduction}
\label{intr}
%%%%%%%%%%%%%%%%%%%%%%%%%%%%%%%%%%%%%%%%%%%%%%%%%%%%%%%%%%%%%%%%%%%%%%%%%%%%%%%%%%%%

Our main goal is the study of the distinguished dynamical role of the
"quantum stationary states" for a model $U(1)$-invariant nonlinear discrete
Schr\"odinger equation
\be\la{S}
   i\dot\psi(x,t)=-\Delta_L\,\psi(x,t)-\de(x)F(\psi(0,t)),\quad x\in\Z.
\ee
Here $F$ is a continuous function, $\delta(x)=\delta_{0x}$ and
$\Delta_L$ stands for the difference Laplacian in $\Z$, defined by
\be\nonumber
   \Delta_L \psi(x)=\psi(x+1)-2\psi(x)+\psi(x-1), \quad x\in\Z
\ee
for functions $\psi: \Z\to\C$.
Physically, equation  (\re{S}) describes the system of the free discrete
Schr\"odinger equation coupled to an oscillator attached at the point $x=0$:
$F$ is a nonlinear ``oscillator force''.

We identify a complex number $\psi=\psi_1+i\psi_2$ with the real two-dimensional
vector $\Psi=(\psi_1,\psi_2)\in\R^2$ and assume that the vector version
$\bf F$ of the oscillator
force $F$ admits a real-valued potential,
\be\la{P}
  {\bf F}(\Psi)=-\na U(\Psi),\quad\Psi\in\R^2,  ~~~~U\in C^2(\R^2).
\ee
Then (\re{S}) is  a Hamiltonian system with Hamiltonian
\be\la{H}
 {\cal H}(\Psi)=\fr 12\langle -\Delta_L\Psi,\Psi\rangle+U(\Psi(0))
=\fr 12\langle\nabla_L\Psi,\nabla_L\Psi\rangle+U(\Psi(0)),
\ee
where $\langle\cdot\rangle$ stands for the scalar product in $l^2(\Z)$,
and $\nabla_L\psi(x)=\psi(x+1)-\psi(x)$.

We assume that
$U(\psi)=u(|\psi|^2)$ with $u\in C^2(\R)$.
Therefore,
by (\re{P}),
\be\la{I}
   F(\psi)= a(|\psi|^2)\psi,\quad\psi\in\C\5,~~~~~~~~a\in C^1(\R),
\ee
where $a(|\psi|^2)$ is real.
Then
$F(e^{i\theta}\psi)=e^{i\theta} F(\psi),\quad\theta\in[0,2\pi] $ and $F(0)=0$
for continuous $F$.
Hence,  $e^{i\theta}\psi(x,t)$ is a solution to (\re{S})  if $\psi(x,t)$ is.
Therefore,  equation (\ref{S}) is $U(1)$-invariant in
the sense of \ci{GSS}, and the N\"other theorem implies the  conservation of $l^2$ norm:
$$
  \Vert\psi(t)\Vert= \Vert\psi(0)\Vert.
$$
Here and below $\Vert\cdot\Vert=\Vert\cdot\Vert_{l^2}$.

The main subject of this paper is an analysis of the special role
played by {\it solitary waves} or {\it nonlinear eigenfunctions},
which are  finite energy solutions of the form
\be\la{SW}
  \psi(x,t)=\psi_\om(x)e^{i\om t},\quad \om\in\R.
\ee
The frequency  $\om$ and the amplitude $\psi_\om(x)$ solve the following
{\it nonlinear eigenvalue problem}:
\be\la{NEP}
  -\om\psi_\om(x)= -\Delta_L\psi_\om(x)-\de(x)F(\psi_\om(0)),\quad x\in\Z
\ee
which follows directly from  (\re{S}) and (\re{I}) since $\om\in\R$.
The  solitary waves constitute a two-dimensional {\it solitary manifold} in
the Hilbert phase space $l^2(\Z)$.

We prove the asymptotics of type
\be\la{sol-as-i}
     \psi(\cdot,t)\sim\psi_{\om_{\pm}}e^{i\om_{\pm}t}+W(t)\Phi_{\pm},
     \quad t\to\pm\infty,
\ee
where $W(t)$ is the dynamical group of the free Schr\"odinger equation,
$\Phi_{\pm}\in l^2(\Z)$ are the corresponding asymptotic scattering states,
and the remainder converges to zero as $\cO(|t|^{-1/2})$ in global norm of
$l^2(\Z)$. The asymptotics hold for the solutions with initial states close to the
{\it stable part} of the solitary manifold, extending the results of
\ci{BKKS,BP,BS,MW96,PW92,PW94} to the equation (\re{S}).

For the first time, the asymptotics of type (\re{sol-as-i}) were established
by Soffer and Weinstein \ci{SW1,SW2} (see also \ci{PW}) for nonlinear
$U(1)$-invariant continuous Schr\"odinger equation with small initial states if the
nonlinear coupling constant is sufficiently small.
The next result was obtained by Buslaev and Perelman \ci{BP}
who proved that the solitary manifold attracts finite energy solutions
of a 1D nonlinear  $U(1)$-invariant translation invariant Schr\"odinger equation
with initial states sufficiently close to the {\it stable part} of the solitary manifold.
For a more lengthy discussion of our motivation, and of previous results in the literature
(\cite{BS,PW92,PW94,PW,SW1,SW2,SW04}) we refer the reader to the introduction of \cite{BKKS}.
The asymptotics of type (\re{sol-as-i}) for nonlinear discrete Schr\"odinger equation
are obtained for the first time in the present paper.

Let us note that we impose conditions which are more general
than the standard ones in the following respects:
we do not hypothesize any spectral properties of the linearized equation,
and do not require any smallness condition on the initial state (only
closeness to the solitary manifold).
This progress is possible on account of the simplicity of our model
which allows an exact analysis of  spectral properties of the linearization.

Let us comment on the general  strategy of our proofs.
We develop the approach \ci{BKKS,BP,IKV05,KK} for our problem. Firstly, we apply the
{\it symplectic projection} onto the solitary manifold to separate the motion along
the solitary manifold and in transversal direction. Secondly, we derive the
{\it modulation equations} for the parameters of the symplectic
projection, and linearize the transversal dynamics at the projection of
the trajectory. The linearized equation is nonautonomous, and this is one of the
fundamental difficulties in the proof. This difficulty is handled by the introduction
of an autonomous  equation (by freezing the time) with an application of the modulation
equations to estimate the resulting additional error terms. A principal role in the rest
of the proof is played by the uniform decay of the frozen linearized dynamics projected
onto the continuous spectrum, and the method of majorants.

The paper is organized as follows.
In \S\ref{swsec} we describe all nonzero solitary waves
and formulate the main theorem.
In \S\ref{subspace} we enumerate some properties of the linearized equation.
In \S\ref{rpsec} we establish the time decay for the linearized
equation in the continuous spectrum.
In \S\ref{modsec} the modulation equations for the parameters of
the soliton are displayed. The decay of the transverse component is proved in
\S\ref{sassec} and \S\ref{prsec}.
In \S\ref{solas-sec} we obtain the soliton asymptotics (\ref{sol-as-i}).
In Appendix we study the resolvent of linearized equation.
%%%%%%%%%%%%%%%%%%%%%%%%%%%%%%%%%%%%%%%%%%%%%%%%%%%%%%%%%%%%%%%%%%%%%%%%%%%%%%%%%%%%%%%%%%%
%%%%%%%%%%%%%%%%%%%%%%%%%%%%%%%%%%%%%%%%%%%%%%%%%%%%%%%%%%%%%%%%%%%%%%%%%%%%%%%%%%%%%%%%%%%%%
\setcounter{equation}{0}
\section{Solitary waves and statement of the main theorem}
\label{swsec}
%%%%%%%%%%%%%%%%%%%%%%%%%%%%%%%%%%%%%%%%%%%%%%%%%%%%%%%%%%%%%%%%%%%%%%%%%%%%%%%%%%%%%%%%%%%%

Our main results describe the large time behavior of the global solutions
whose existence is guaranteed by the following theorem
\begin{theorem}\label{locex}\cite[Theorem 3.1]{K08}\\
i)
   Let conditions (\re{P}) and  (\re{I}) hold. Then for any $\psi_0\in l^2=l^2(\Z)$
   there exist a unique solution $\psi\in C_b(\R,l^2)$
   to the equation (\ref{S}) with initial condition $\psi(x,0)=\psi_0(x)$.\\
ii)
   The value of energy functional and the norm of solution is conserved:
   \be\la{const}
   {\cal H}(\psi(t))={\cal H}(\psi_0),\quad \Vert\psi(t)\Vert=\Vert\psi_0\Vert,
   \quad t\in\R.
   \ee
\end{theorem}
%%%%%%%%%%%%%%%%%%%%%%%%%%%%%%%%%%%%%%%%%%%%%%%%%%%%%%%%%%%%%%%%%%%%%%%%%%%%%%%%%%%%%%%%%%
In \cite{K08} we give a complete analysis of the set of all solitary waves.
There exist two different sets of nonzero solitary waves.
The first set $\cS_+$  corresponds to $\om\in(0,\infty)$
and the second set $\cS_-$  corresponds  to $\om\in (-\infty,-4)$.

Denote by $k(\om)$ the positive solution of the equation $\cosh k=|\om+2|/2$.
\begin{lemma}\cite[Lemma 4.1]{K08})
The sets of all nonzero solitary waves is given by
$$
\cS_+=\Bigl\{\psi_{\om}e^{i\theta}\!=\!C e^{i\theta-k(\om) |x|}:
\;\om\in(0,\infty),\;  C>0,\;
\sinh k(\om) =a(C^2)/2>0,\;\theta\in[0,2\pi]\Bigr\}
$$
$$
\cS_-\!=\!\Bigl\{\psi_{\om}e^{i\theta}\!=C(\!-1\!)^{|x|} e^{i\theta-k(\om) |x|}\!:
\;\om\in (-\infty,-4),\;  C\!>\!0,\;
\sinh k(\om)\!=\!-a(C^2)/2>\!0,\;\theta\in \![0,2\pi]\Bigr\}
$$
\end{lemma}
%%%%%%%%%%%%%%%%%%%%%%%%%%%%%%%%%%%%%%%%%%%%%%%%%%%%%%%%%%%%%%%%%%%%%%%%%%%
\bc
The set $\cS_+$ resp. $\cS_-$ is a smooth manifold with the co-ordinates
$\theta\in\R\mod 2\pi$
and $C>0$ such that $a(C^2)>0$ resp. $a(C^2)<0$.
\ec
%%%%%%%%%%%%%%%%%%%%%%%%%%%%%%%%%%%%%%%%%%%%%%%%%%%%%%%%%%%%%%%%%%%%%%%%%%%%%%%%
In the case of polynomial $F$ the condition on $C$ means that $C$ is restricted
to lie in a set which is a finite union of one-dimensional intervals.
The value $C=0$ corresponds to the zero function $\psi_\om(x)=0$ which
is always a solitary  wave since $F(0)=0$, and for $\om\in [-4,0]$
only the zero solitary wave exists.

In (\ci {K08}) we proved that the parameters $\theta, \om$ locally also are smooth
coordinates on $\cS_{\pm}$ at the points with $a'(C^2)\ne 0$.
We will also need the following result from \cite{K08}.
\begin{lemma}\la{int-dif}
For $C> 0$, $a>0$, and $a'\not=(4a+a^3)/(4C^2)$ we have
\be\la{syc2}
  \pa_\om\int|\psi_\om(x)|^2 dx\not=0.
\ee
\end{lemma}
%%%%%%%%%%%%%%%%%%%%%%%%%%%%%%%%%%%%%%%%%%%%%%%%%%%%%%%%%%%%%%%%%%%%%%%%%%%%%%%%%%%
The soliton solution  is a trajectory $\psi_{\om(t)}(x)e^{i\theta(t)}$,
where the parameters  satisfy the equation $\dot\theta=\om$, $\dot\om=0$.
The solitary waves $e^{i\theta}\psi_\om(x)$ map out in time an orbit
of the $U(1)$ symmetry group. This group acts on the phase space $l^2({\Z})$ preserving
the Hamiltonian ${\cal H}$.

The real form of the solitary wave is $e^{j\theta}\Phi_\om$ where $\Phi_\om=(\psi_{\omega}(x),0)$,
and $j$ is the $2\times 2$ matrix
\be\la{j}
   j=\left(
   \ba{rr}
   0  &  -1\\
   1  &   0
   \ea
   \right)~~~~~~~~~~~~~~~~~~~~~~~~~~~~~~~~
\ee
Linearization at the solitary wave $e^{j\theta}\Phi_\om$ leads to the operator
(cf. \ci{BKKS, BS})
\be\la{B}
{\bf B}=-\Delta_L+\om-\de(x)[a(C^2)+2a'(C^2)C^2P_1]=
\left(
   \ba{cc}
   {\bf D}_1  &       0    \\
        0     &   {\bf D}_2
   \ea
   \right),
\ee
where $P_1$ is the projector in $\R^2$ acting as
$\left(\ba{l}\chi_1\\ \chi_2\ea\right)\mapsto \left(\ba{l}\chi_1\\ 0\ea\right)$,
$$
   {\bf D}_1=-\Delta_L+\om-\de(x)[a+2a'C^2],\quad
   {\bf D}_2=-\Delta_L+\om-\de(x)a.
$$
Let $\bC=j^{-1}\bB$.
We will show in Appendix A that the continuous spectrum of $\bC$  coincides with
$\cC_-\cup\cC_+=[-i(\om+4),-i\om]\cup [i\om, i(\om+4)]$.
The point 0 belongs to the discrete spectrum, and the dimension of its invariant subspace
is at least 2. If $a'\not=(4a+a^3)/(2C^2)$ then the invariant subspace associated to
the eigenvalue $\lam=0$ is of dimension exactly 2. We will analyze only the solitary waves
with $a'\not\in\{0;(4a+a^3)/(2C^2)\}$. We assume also more specific condition
\bd\la{spcon}
  We say the solitary wave $\psi_{\om}(x)e^{i\theta}$,
  satisfies the spectral condition $\{SP\}$ if\\
  1) $a'\not\in\{0;(4a+a^3)/(2C^2)\}$\\
  2) there is no  eigenvalue except $\lam=0$.
\ed
In Appendix A we give an example  when the condition $\{SP\}$ holds.
If condition $\{SP\}$ is true for a fixed value $\om_0$, it is also true for
values $\om$ in a small interval centered at $\om_0$. The condition $\{SP\}$
ensures orbital stability of solitary waves.
%%%%%%%%%%%%%%%%%%%%%%%%%%%%%%%%%%%%%%%%%%%%%%%%%%%%%%%%%%%%%%%%%%%%%%%%%%%%%%%%%%%%%%%%%%

The functional spaces we are going to consider are the weighted Banach spaces
$l^p_{\beta}=l^p_{\beta}(\Z)$, $p\in [1,\infty)$, $\beta\in\R$ of complex valued
functions with the norm
\be\la{norm}
  \Vert u\Vert_{l^p_{\beta}}=\Vert (1+|x|)^{\beta} u(x)\Vert_{l^p}.
\ee
Let us denote by $W(t)$ the dynamical group of the free Schr\"odinger equation.
Our main  theorem is the following:
%%%%%%%%%%%%%%%%%%%%%%%%%%%%%%%%%%%%%%%%%%%%%%%%%%%%%%%%%%%%%%%%%%%%%%%%%%%%%%%%%%%%%%%%%%%%
\begin{theorem}\label{main}
   Let conditions (\re{P}),  and  (\re{I}) hold,
   $\beta\ge 2$ and $\psi(x,t)\in C(\R,l^2)$ be the solution to the equation (\ref{S})
   with initial value $\psi_0(x)=\psi(x,0)\in l^2\cap l^1_\beta$ which
   is close to a solitary wave
   $\psi_{\om_0}e^{i\theta_0}$
   \be\la{close}
     d:=\Vert\psi_0-\psi_{\om_0}e^{i\theta_0}\Vert_{l^2\cap l^1_{\beta}}\ll 1.
   \ee
   Assume further that the spectral condition $\{SP\}$ holds for the solitary wave with $\om=\om_0$.
   Then for $d$ sufficiently small the solution admits the following asymptotics:
   \be\la{sol-as}
     \psi(\cdot,t) = \psi_{\om_{\pm}}e^{i\om_{\pm}t}+W(t)\Phi_{\pm}
     +r_{\pm}(t),\quad t\to\pm\infty,
   \ee
   where $\Phi_{\pm}\in l^2(\Z)$ are the corresponding asymptotic
   scattering states, and
   \be\la{rate}
   \Vert r_{\pm}(t)\Vert=\cO(|t|^{-1/2}),\quad t\to\pm\infty.
   \ee
\end{theorem}

%%%%%%%%%%%%%%%%%%%%%%%%%%%%%%%%%%%%%%%%%%%%%%%%%%%%%%%%%%%%%%%%%%%%%%%%%%%%%%%%%%%%%%%%%
\setcounter{equation}{0}
\section{Linearized evolution}
\label{subspace}
%%%%%%%%%%%%%%%%%%%%%%%%%%%%%%%%%%%%%%%%%%%%%%%%%%%%%%%%%%%%%%%%%%%%%%%%%%%%%%%%%%%%%%%%%
In this section we summarize the properties of the linearized evolution which will be needed.
The proof of these properties one can find in Appendix A and in \cite{BKKS}.
The linearized equation reads
\be\la{lin4}
   \dot\chi(x,t)={\bf C}\chi(x,t),~~~~~{\bf C}:=j^{-1}\bB=
   \left(
   \ba{rr}
    0          &    {\bf D}_2\\
   -{\bf D}_1  &        0
   \ea
   \right).
\ee
Theorem \ref{locex} generalizes to the equation (\ref{lin4}):
the equation admits unique solution $\chi(x,t)\in C_b(\R,l^2)$
for every initial function $\chi(x,0)=\chi_0\in l^2$.
Denote by $e^{{\bf C}t}$ the dynamical group of equation (\ref{lin4})
acting in the space $l^2$. Then (\ref{const}) implies that
\be\la{t-small}
  \Vert e^{{\bf C}t}\chi_0\Vert=\Vert\chi_0\Vert,\quad t\in \R.
\ee
The resolvent $\bR(\lam):=(\bC-\lam)^{-1}$ is an  integral operator
with matrix valued integral kernel (see Appendix A)
\be\la{nR}
  {\bf R}(\lam,x,y)=\Gamma(\lam,x,y)+P(\lam,x,y),
\ee
where
\be\la{tG}\!
   \Gamma(\lam,x,y)=\left(\!\! \ba{cc}
   \ds\frac 1{4\sin\theta_+}  &  -\ds\frac 1{4\sin\theta_-}\\
   \ds\frac i{4\sin\theta_+}  &  \ds\frac i{4\sin\theta_-}
   \ea\!\!\right)\!\!\left(\!\! \ba{cc}
   e^{i\theta_+|x-y|}-e^{i\theta_+(|x|+|y|)}  &  -i(e^{i\theta_+|x-y|}-e^{i\theta_+(|x|+|y|)})\\\\
   e^{i\theta_-|x-y|}-e^{i\theta_-(|x|+|y|)}  &   i(e^{i\theta_-|x-y|}-e^{i\theta_-(|x|+|y|)})
   \ea\!\!\right)
\ee
\be\la{tP}\!
   P(\lam,x,y)=\frac 1{2D} \left(\!\! \ba{cc}
   e^{i\theta_+|x|}    &   e^{i\theta_-|x|}\\
   ie^{i\theta_+|x|}   &  -ie^{i\theta_-|x|}
   \ea\!\!\right)\!\!
   \left(\!\! \ba{cc}
   i\al-2\sin\theta_-       &    i\beta\\
   -i\beta         &   -i\al+2\sin\theta_+
   \ea\!\!\right)\!\!
   \left(\!\! \ba{cc}
   e^{i\theta_+|y|}      &  -ie^{i\theta_+|y|}      \\
   e^{i\theta_-|y|}       &  ie^{i\theta_-|y|}
   \ea\!\!\right)
\ee
Here $\theta_\pm(\lam)$ is
the  root of $2\cos\theta_{\pm}=\om+2\pm i\lam$ defined with cuts in the complex  $\lam$ plane
so that $\theta_\pm(\lambda)$ is analytic on ${\C}\setminus {\cal C}_\pm$, and
${\rm Im}\5 \theta_\pm(\lambda)>0 $ for $\lam\in\C\setminus {\cal C}_\pm$.
The constants $\al$, $\beta$ and the determinant $D=D(\lam)$ are given by the formulas
$$
 \al=a+a'C^2,\; \beta=a'C^2,\; D=2i\al(\sin\theta_++\sin\theta_-)
 -4\sin\theta_+\sin\theta_-+\al^2-\beta^2.
$$
%%%%%%%%%%%%%%%%%%%%%%%%%%%%%%%%%%%%%%%%%%%%%%%%%%%%%%%%%%%%%%%%%%%%%%%%%%%%%%%%%%%%%%%%%
The poles of the resolvent correspond to the roots of the determinant $D(\lam)$.\\
If spectral condition $\{SP\}$ holds
then the determinant has the only root $\lam =0$ with the multiplicity $2$.

Observe that
(\ref{NEP}) and its derivative in $\om$ give the following identities:
\be\la{diffeq}
  {\bf D}_2\psi_\omega=0\qquad {\bf D}_1(\partial_\omega\psi_\omega)=-\psi_\omega.
\ee
These formulae imply that the vectors $j\Phi_\om$ and $\pa_\om\Phi_\om$ lie in the generalized
two dimensional null space $X^0$ of the non-self-adjoint operator ${\bf C}$ defined in
(\ref{lin4}) and
\be\la{CT}
   {\bf C} j\Phi_\om=0\qquad {\bf C}\pa_\om\Phi_\om=j\Phi_\om.
\ee
The symplectic form $\Omega$ for the real vectors $\psi$ and $\eta$ is defined by
\be\la{symp}
   \Om(\psi,\eta)=\langle\psi,j\eta\rangle=\sum\limits_{\Z}(\psi_1\eta_2-\psi_2\eta_1).
\ee
By Lemma \ref{int-dif}
\be\label{nondeg}
   \Omega(j\Phi_\om,\pa_\om\Phi_\om)=-\frac{1}{2}\pa_\om\int |\psi_\om|^2dx\not=0.
\ee
Hence, the symplectic form $\Omega$ is nondegenerate on $X^0$,
i.e. $X^0$ is a symplectic subspace. Therefore, there exists a symplectic projection
operator ${\bf P}^0$ from $l^2$ onto $X^0$ represented by the formula
\be\label{defsp}
{\bf P}^0\psi=\frac 1{\langle\Phi_\om,\partial_\om\Phi_\om\rangle}
  [\langle\psi,j\partial_\om\Phi_\om\rangle j\Phi_\om
  +\langle \psi,\Phi_\om\rangle\partial_\om\Phi_\om]
\ee
Denote by $\bP^c=1-\bP^0$  the symplectic projector onto
the continuous spectral subspace.
\br\la{red}
On the generalized null space itself ${\bf C}^2=0$ by (\ref{CT}), and so the
semigroup $e^{t{\bf C}}$ reduces to $1+{\bf C}t$ as usual for the exponential of the
nilpotent part of an operator.
\er
%%%%%%%%%%%%%%%%%%%%%%%%%%%%%%%%%%%%%%%%%%%%%%%%%%%%%%%%%%%%%%%%%%%%%%%%%%%%%%%%%%%%%%%%%%
%%%%%%%%%%%%%%%%%%%%%%%%%%%%%%%%%%%%%%%%%%%%%%%%%%%%%%%%%%%%%%%%%%%%%%%%%%%%%%%%%%%%%%%%
\setcounter{equation}{0}
\section{Time decay in continuous spectrum}
\label{rpsec}
%%%%%%%%%%%%%%%%%%%%%%%%%%%%%%%%%%%%%%%%%%%%%%%%%%%%%%%%%%%%%%%%%%%%%%%%%%%%%%%%%%%%%%%%%
Due to  Remark \ref{red} we see that the solutions
$\chi(t)= e^{{\bf C}t}\chi_0$ of the linearized equation (\ref{lin4})
do not decay as $t\to\infty$ if $\bP^0\chi_0\ne 0$.
On the other hand, we do expect time decay of ${\bf P}^c\chi(t)$,
as a consequence of the Laplace representation  for
${\bf P}^c e^{{\bf C}t}$:
\be\la{Pc-rep}
 {\bf P}^c e^{{\bf C}t}
  =-\frac 1{2\pi i}\int\limits_{{\cal C}_+\cup{\cal C}_-}
  e^{\lam t} ({\bf R}\bigl(\lam+0)-{\bf R}(\lam-0)\bigr)~d\lam.
\ee
The decay for the oscillatory integral is obtained from the analytic
properties of $\bR(\lam)$ for $\lam\in{\cal C}_+\cup{\cal C}_-$.

In Appendix A we will show  that $D(\lam)\ne 0$ for $\lam\in{\cal C}_+\cup{\cal C}_-$.
Clearly in order to understand the decay of ${\bf P}^c e^{t{\bf C}}$,
it is crucial to study the behavior of $ {\bf R}(\lam,x,y)$ near the
branch points $\lam=\pm i\om$ and $\lam=\pm i(\om+4)$ (where $\sin\theta_\pm$ vanish).

We deduce time decay for the group ${\bf P}^c e^{t{\bf C}}$
by means of the following  version of Lemma 10.2 from \cite{JK},
which is itself based on Zygmund's lemma \ci[p.45]{Z}.

Let ${\cal F}:[a,b]\to {\bf B}$ be a $C^2$ function with values in a
Banach space ${\bf B}$.
Let us define the ${\bf B}$-valued function
$$
I(t)=\int\limits^a_b e^{-it\nu}{\cal F}(\nu)~d\nu.
$$
\begin{lemma}\la{K}
 Suppose that ${\cal F}(a)={\cal F}(b)=0$, and
 ${\cal F}''\in L^1(a+\delta,b-\delta; {\bf B})$ for some $\delta>0$.
 Moreover,
 $$
    {\cal F} ''(a+\zeta),\;{\cal F} ''(b-\zeta)={\cal O}(\zeta^{p-2}),\quad \zeta\downarrow 0
 $$
  in the norm of $~{\bf B}$ for some $p\in (0,1)$.
  Then $I(t)\in C_b(a+\ve,b-\ve;\bB)$ for any $\ve>0$, and
  $$
  I(t) ={\cal O}(t^{-1-p})
  \quad{\rm as}\quad t\to\infty\quad {\rm in~ the~ norm~ of}\quad\bB.
  $$
\end{lemma}
We will apply Lemma \ref{K} to the function
${\cal F}(\lam)={\bf R}(\lam+0)-{\bf R}(\lam-0)$ with values in the Banach space
${\cal B} =B(l^1_\beta,l^{\infty}_{-\beta})$ ,
the space of continuous linear maps $l^1_\beta\to l^\infty_{-\beta}$ for any $\beta\ge 2$.
\begin{theorem}\la{TD}
  Assume that the spectral condition $\{SP\}$ holds so that
  $\lambda=0$ is the only point in the discrete spectrum of the
  operator $\bC=\bC(\om)$. Then
for $\beta\ge 2$
  \be\la{decayR}
    \Vert {\bf P}^c e^{{\bf C}t}\Vert_{\cal B}
    ={\cal O}(t^{-3/2}),\quad t\to\infty.
  \ee
\end{theorem}
First we
use the formulas (\ref{Pc-rep}) and (\ref{nR}) to obtain
\be\la{Z}
  -2\pi i{\bf P}^c e^{{\bf C}t}=
  \int\limits_{{\cal C}_+\cup{\cal C}_-}\!
  e^{\lam t}(\Gamma(\lam+0)-\Gamma(\lam-0))\5d \lam
 ~~ +\int\limits_{{\cal C}_+\cup{\cal C}_-}\!
  e^{\lam t}(P(\lam+0)-P(\lam-0))\5 d\lam
\ee
Next we apply Lemma \ref{K} to each summand in the RHS of (\ref{Z})
separately. Then Theorem \ref{TD} immediately follows from the  two lemmas below.
%%%%%%%%%%%%%%%%%%%%%%%%%%%%%%%%%%%%%%%%%%%%%%%%%%%%%%%%%%%%%%%%%%%%%%%%%%%%%%%%%%%%%%%%%%%%%%%%
\begin{lemma}\la{decayG}
If the assumption of Theorem \ref{TD} hold then
  \be
    \int\limits_{{\cal C}_+\cup{\cal C}_-}
    e^{\lam t}(\Gamma(\lam+0)-\Gamma(\lam-0))~d\lam
    ={\cal O}(t^{-3/2}),\quad t\to\infty
  \ee
  in the norm $\cal B$.
\end{lemma}
{\bf Proof }
We consider only the integral over ${\cal C}_+$ since the integral over ${\cal C}_-$
can be handled in the same way.
The points $\lam=i\om$ and $\lam=i(\om+4)$ are the branch points for $\theta_+$, therefore,
if $\lam\in{\cal C}_+$ then since $\theta_-$ is continuous across ${\cal C}_+$
$$\Gamma(\lam+0)-\Gamma(\lam-0)=\Gamma^+(\lam+0)-\Gamma^+(\lam-0),$$
where $\Gamma^+$ is the  sum of those terms in $\Gamma$ which involve $\theta_+$.
Let us consider, for example,  $\Gamma^+_{11}$.
The expression (\re{tG}) implies for $y>0$ that
$$
  \Gamma^+_{11}(\lam,x,y)= \left\{ \ba{ll}
  0,  \!&\!x\le0,\\\\
  \ds\frac{e^{i\theta_+ y}(e^{-i\theta_+ x}-e^{i\theta_+ x})}{4\sin\theta_+},  \!&\!0\le x\le y,\\\\
  \ds\frac{e^{i\theta_+ x}(e^{-i\theta_+ y}-e^{i\theta_+ y})}{4\sin\theta_+},  \!&\!x\ge y.
  \ea \right.
$$
For $\lam\in {\cal C}_+$, the root
$\theta_+$ is real, and $\theta_+(\lam+0)=-\theta_+(\lam-0)$.
Then, for $y>0$,
\beqn
  \Gamma^+_{11}(\lam+0,x,y)-\Gamma^+_{11}(\lam-0,x,y)&=&
  -\Theta(x)\ds\frac{\sin\theta_+|x|\sin\theta_+|y|}{\sin\theta_+},
\eeqn
where  $\Theta(x)=1$ for $x>0$ and zero otherwise. Recall that
$$
\theta_+=\theta_+(\nu)=\arccos\fr{\om+2-\nu}{2},
\quad\sin\theta_+=\fr 12\sqrt{(-\om+\nu)(4+\om-\nu)},\quad \nu=-i\lam.
$$
Let us calculate the second derivative of the function
$f(\nu)=\ds\frac{\sin\theta_+|x|\sin\theta_+|y|}{\sin\theta_+}$:
\beqn\nonumber
f''(\nu)\!\!\!&=&\!\!\!\fr{-2(|x|^2+|y|^2-1)\sin\theta_+|x|\sin\theta_+|y|
+4|x||y|\cos\theta_+|x|\cos\theta_+|y|}{(-\om+\nu)^{3/2}(4+\om-\nu)^{3/2}}\\
\nonumber
\!\!\!&-&\!\!\!\fr{8(|x|\cos\theta_+|x|\sin\theta_+|y|
+|y|\sin\theta_+|x|\cos\theta_+|y|)\cos\theta_+}{(-\om+\nu)^2(4+\om-\nu)^2}
+\fr{16\sin\theta_+|x|\sin\theta_+|y|\cos^2\theta_+}{(-\om+\nu)^{5/2}(4+\om-\nu)^{5/2}}\\
\nonumber
\!\!\!&+&\!\!\!\Bigg[\frac{|x|\cos\theta_+|x|\sin\theta_+|y|+|y|\sin\theta_+|x|\cos\theta_+|y|}
{(-\om+\nu)^2(4+\om-\nu)^2}
-\fr{\sin\theta_+|x|\sin\theta_+|y|\cos\theta_+}{(-\om+\nu)^{5/2}(4+\om-\!\nu)^{5/2}}\Bigg]\!
(\om+2-\!\nu)
\eeqn
Since
$$
\sin\theta_+|z|\sim |z|\sqrt{\nu-\om},\; \nu\to\om\quad{\rm and}\quad
\sin\theta_+|z|\sim |z|\sqrt{4+\om-\nu},\; \nu\to\om+4,
$$
then
$$
|f''(\nu)|\le\ds\frac{C(1+|x|^2)(1+|y|^2)}{(\nu-\om)^{3/2}},\;\om<\nu<\om+1,
$$
and
$$
|f''(\nu)|\le\ds\frac{C(1+|x|^2)(1+|y|^2)}{(4+\om-\nu)^{3/2}},\;\om+3<\nu<\om+4,
$$
For  $y<0$ an identical calculation leads to the same bound.
Therefore the operator valued function ${\cal F}(\nu)=\Gamma^+_{11}(\lam+0)-
\Gamma^+_{11}(\lam-0)$
satisfies the conditions  of Lemma \re{K} with
$a=\om$, $b=4+\om$, $p=1/2$ and ${\bf B}={\cal B}$.~~\bo
\\
%%%%%%%%%%%%%%%%%%%%%%%%%%%%%%%%%%%%%%%%%%%%%%%%%%%%%%%%%%%%%%%%%%%%%%%%%%%%%%%%%%%%%%%%%%%%%%
Next we consider the second summand in the RHS of (\ref{Z}).
\begin{lemma}\la{decayP}
   In the situation of Theorem \ref{TD}
   \be
      \int\limits_{{\cal C}_+\cup{\cal C}_-}
      e^{\lam t}(P(\lam+0)-P(\lam-0))~d\lam={\cal O}(t^{-3/2}),
   \ee
   in the norm ${\cal B}$.
\end{lemma}
{\bf Proof }  We consider only the integral over
${\cal C}_+$ and one component of the matrix $P$, for example,
$P_{11}$:
$$
P_{11}(\lam,x,y)=\!\ds\frac{(i\al\!-\!2\sin\theta_-\!)e^{i\theta_+(|x|+|y|)}
+(2\sin\theta_-\!-\!i\al)e^{i\theta_-(|x|+|y|)}+i\beta (e^{i\theta_-|y|+i\theta_+|x|}\!-\!
e^{i\theta_+|y|+i\theta_-|x|})}{2i\al(\sin\theta_++\sin\theta_-)
-4\sin\theta_+\sin\theta_-+\al^2-\beta^2}
$$
Denote $\zeta=-\om-i\lam$, then
$\sin\theta_+=\sqrt{\zeta}\sqrt{4-\zeta},\;\sin\theta_-=\sqrt{-2\om-\zeta}\sqrt{4+2\om+\zeta}$.
The Taylor expansion in $\sqrt\zeta$ as $\zeta\to 0,\quad \rIm\zeta\ge 0$ implies
$$
P_{11}(i\om+i\zeta,x,y)=P_0+P_1(x,y)\zeta^{1/2}+P_2(x,y){\cal O}(\zeta),
$$
where $|P_j(x,y)|\le C_j(1+|x|^j)(1+|y|^j),\;j=1,2.$
Therefore, if $\lam=i\om+i\zeta\in{\cal C}_+$ then
$${\cal F}(\om+\zeta)=P_{11}(i\om+i\zeta+0)-P_{11}(i\om+i\zeta-0)
={\cal O}(\zeta^{1/2}),\quad \zeta\to 0$$
in the norm of ${\cal B}$.
Similarly, differentiating two times the function $P_{11}(\lam,x,y)$
in $\lam$, we obtain that
$$
{\cal F}''(\om+\zeta)=-P''_{11}(i\om+i\zeta+0)+P''_{11}(i\om+i\zeta-0)
={\cal O}(\zeta^{-3/2}),\quad \zeta\to 0
$$
in the norm of ${\cal B}$.
In the same way
$$
{\cal F}''(4+\om-\zeta)=-P''_{11}(i(4+\om-\zeta)+0)+P''_{11}(i(4+\om-\zeta)-0)
={\cal O}(\zeta^{-3/2}),\quad \zeta\to 0
$$
Therefore, the function ${\cal F}(\nu)$, $\nu=-i\lam$
satisfies the conditions  of Lemma \re{K} with $p=1/2$
and ${\bf B}={\cal B}$.~~\bo
%%%%%%%%%%%%%%%%%%%%%%%%%%%%%%%%%%%%%%%%%%%%%%%%%%%%%%%%%%%%%%%%%%%%%%%%%%%%%%%%%%%%

%%%%%%%%%%%%%%%%%%%%%%%%%%%%%%%%%%%%%%%%%%%%%%%%%%%%%%%%%%%%%%%%%%%%%%%%%%%%%%%%%%%%%%%%%%%%%
%%%%%%%%%%%%%%%%%%%%%%%%%%%%%%%%%%%%%%%%%%%%%%%%%%%%%%%%%%%%%%%%%%%%%%%%%%%%%%%%%%%%
\setcounter{equation}{0}
\section{Modulation equations }
\label{modsec}
%%%%%%%%%%%%%%%%%%%%%%%%%%%%%%%%%%%%%%%%%%%%%%%%%%%%%%%%%%%%%%%%%%%%%%%%%%%%%%%%%%%%%%

In this section we present the modulation equations which allow a
construction of solutions $\psi(x,t)$ of  equation (\re{S})
close at each time $t$ to a soliton i.e. to one of the functions
$\psi_\om(x)$ in the set $\cS_+\cup\cS_-$ described in section \ref{swsec}
with time varying (``modulating'') parameters $(\omega,\theta)=(\omega(t),\theta(t))$.
Let us rewrite (\re{S}) in the  real form
\be\la{SV}
  j\dot\psi(x,t)= -\psi''(x,t)-\de(x)\bF(\psi(0,t)),
\ee
as an equation for $\psi(x,t)\in\R^2$
with $\bF(\psi)\in\R^2$ which is the real vector version of $F(\psi)\in\C$.
Then $\psi(x,t)=e^{j\theta(t)}\Phi_{\om(t)}(x)$ is a solution of (\ref{SV}) if and
only if $\dot\theta=\om$ and $\dot\om=0$.

We look for a solution to (\ref{SV}) in the form
\be\la{sol}
  \psi(x,t)=e^{j\theta(t)}\bigl(\Phi_{\om(t)}(x)+\chi(x,t)\bigr)=
  e^{j\theta(t)}\Psi(x,t),\quad \Psi(x,t)=\Phi_{\om(t)}(x)+\chi(x,t).
\ee
Since this is a solution of (\ref{SV}) as long as
$\chi\equiv 0$ and $\dot\theta=\om$ and $\dot \omega=0$
it is natural to look for solutions  in which $\chi$ is small and
$$\theta(t)=\int_0^t\om(s)ds+\gamma(t)$$ with $\gamma$ treated perturbatively.
Observe that so far this representation is underdetermined since for any
$\bigl(\om(t),\theta(t)\bigr)$ it just amounts to a
definition of $\chi$; it is made unique by restricting
$\chi(t)$ to lie in the image of the projection operator
onto the continuous spectrum ${\bf P}^c_t={\bf P}^c(\omega(t))$ or equivalently that
\be\la{req}
  \;{\bf P}^0_t \chi(t)=0,\;\;{\bf P}^0_t={\bf P}^0(\omega(t))=I-{\bf P}^c(\omega(t))
\ee
Now we give a system of {\it modulation equations} for
$\omega(t),\;\gamma(t)$ which ensure the conditions (\ref{req})
are preserved by the time evolution.
%%%%%%%%%%%%%%%%%%%%%%%%%%%%%%%%%%%%%%%%%%%%%%%%%%%%%%%%%%%%%%%%%%%%%%%%%%%%%%%%%%%%%%%%%%
\begin{lemma}\label{BS}
  (i) Assume given a solution of (\ref{SV}) with regularity
  as described in Theorem \ref{locex}, which can be
  written in the form (\ref{sol}) -(\ref{req}) with continuously differentiable
  $\omega(t),\;\theta(t)$. Then
  \be\la{SV2}
   \dot\chi={\bf C}\chi-\dot\om\pa_\om\Phi_\om+\dot\gamma j^{-1}(\Phi_\om+\chi)+{\bf Q}
  \ee
  where
  ${\bf Q}(\chi,\om)=-\delta(x) j^{-1}\bigl({\bF}(\Phi_\om+\chi)
  -{ \bF}(\Phi_\om)-{ \bF}'(\Phi_\om)\chi\bigr)$, and
  \beqn\la{omega}
    \dot\omega&=&\frac{\langle {\bf P}^0{\bf Q},\Psi\rangle}
    {\langle\pa_\om\Phi_{\om}-\pa_\om {\bf P}^0\chi,\Psi\rangle}\\
    \dot\gamma&=&\frac{\langle j{\bf P}^0(\pa_\om\Phi_{\om}-\pa_\om {\bf P}^0\chi),
    {\bf P}^0{\bf Q}\rangle}
    {\langle\pa_\om \Phi_\om-\pa_\om {\bf P}^0\chi,\Psi\rangle},\la{gamma}.
  \eeqn
  where ${\bf P}^0={\bf P}^0(\om(t))$ is the projection operator
  defined in (\ref{defsp}) and $\pa_\om{\bf P}^0=\pa_\om{\bf
  P}^0(\om)$ evaluated at $\om=\om(t)$.

  (ii) Conversely given $\psi$ a solution of (\ref{SV}) as in Theorem \ref{locex} and
  continuously differentiable functions $\omega(t),\;\theta(t)$ which satisfy
  (\ref{omega})-(\ref{gamma}) then $\chi$ defined by (\ref{sol}) satisfies (\ref{SV2})
  and the condition (\ref{req}) holds at all times if it holds initially.
\end{lemma}
\Pr
This can be proved as in \ci[Prop.2.2]{BS}.\bo

It remains to show, for appropriate initial data close to a soliton, that there exist
solutions to (\ref{omega})-(\ref{gamma}), at least locally. To achieve this observe
that if the spectral condition $\{SP\}$ holds then by Lemma \ref{int-dif}
the denominator appearing on the right hand side of (\ref{omega}) and (\ref{gamma})
does not vanish for  small $\Vert\chi\Vert_{l^{1}_{\beta}}$. This is because
\be\la{int-dif1}
  \langle \pa_\om\psi_\om,\psi_\om\rangle=\frac{1}{2}\partial_\om\int|\psi_\om|^2dx\not=0
\ee
as discussed in section \ref{swsec}.
This has the consequence that the orthogonality conditions really can be satisfied for
small $\chi$ because they are equivalent to a locally well posed set of ordinary
differential equations for $t\to(\theta(t),\omega(t))$.
This implies the following corollary:
%%%%%%%%%%%%%%%%%%%%%%%%%%%%%%%%%%%%%%%%%%%%%%%%%%%%%%%%%%%%%%%%%%%%%%%%%%%%%%%%%%%%%%%%
\begin{cor}\label{ME}
(i) In the situation of (i) in the previous lemma assume that
(\ref{int-dif1}) holds.
If $\Vert\chi\Vert_{l^p_\beta}$ is sufficiently small for some
$p,\beta$  the right hand sides of
(\ref{omega}) and (\ref{gamma}) are smooth in $\theta,\om$ and there exists
continuous ${\cal R}={\cal R}(\omega, \chi)$ such that
$$
  |\dot\gamma(t)|\le{\cal R}|\chi(0,t)|^2,\qquad
  |\dot\omega(t)|\le{\cal R}|\chi(0,t)|^2.
$$

(ii) Assume given $\psi$, a solution of (\ref{SV}) as in Theorem \ref{locex}.
If $\om_0$ satisfies (\ref{int-dif1}) and
$\chi(x,0)=e^{-j\theta_0}\psi(x,0)-\Phi_{\om_0}(x)$
is small in  some $l^p_\beta$ norm and satisfies (\ref{req}) there is a time interval
on which there exist $C^1$ functions $t\mapsto\bigl(\om(t),\gamma(t)\bigr)$ which
satisfy (\ref{omega})-(\ref{gamma}).
\end{cor}

%%%%%%%%%%%%%%%%%%%%%%%%%%%%%%%%%%%%%%%%%%%%%%%%%%%%%%%%%%%%%%%%%%%%%%%%%%%%%%%%%%%%
%%%%%%%%%%%%%%%%%%%%%%%%%%%%%%%%%%%%%%%%%%%%%%%%%%%%%%%%%%%%%%%%%%%%%%%%%%%%%%%%%%%%
\section{Time decay for the transversal dynamics}
\setcounter{equation}{0}
\label{sassec}
%%%%%%%%%%%%%%%%%%%%%%%%%%%%%%%%%%%%%%%%%%%%%%%%%%%%%%%%%%%%%%%%%%%%%%%%%%%%%%%%%%%%%
Let us represent the initial data $\psi_0$ in a convenient form for
application of the modulation equations: the next Lemma will allow us
to assume that (\ref{req}) holds initially without loss of generality.
\begin{lemma}\la{in-cond}
In the situation of Theorem \ref{main} there exists a solitary wave
$\psi_{{\tilde\om}_0}$
satisfying the spectral condition $\{SP\}$ such that  in vector form
$$
\psi_0(x)=e^{j{\tilde\theta}_0}(\Phi_{{\tilde\om}_0}(x)+\chi_0(x)),
\quad  \Phi_{{\tilde\om}_0}=( \psi_{{\tilde\om}_0},0),
$$
and for $\chi_0(x)$ we have
 \be\la{innit}
   {\bf P}^0({\tilde\omega}_0)(\chi_0)=0,
\ee
and
  $$
    \Vert\chi_0\Vert_{l^1_{\beta}\cap l^2}= \tilde d=O(d)\quad\hbox{as}\; d\to 0.
  $$
\end{lemma}
{\bf Proof}
This can be proved as in \ci[Lemma 10.1]{BKKS}
by a standard application of the implicit function theorem.
\bo

In Section \ref{solas-sec} we will show that our main Theorem \ref{main}
can be derived from the following time decay of the transversal component $\chi(t)$:
%%%%%%%%%%%%%%%%%%%%%%%%%%%%%%%%%%%%%%%%%%%%%%%%%%%%%%%%%%%%%%%%%%%%%%%%%%%%%%%%%%%%%%%%5
\begin{theorem}\la{yest}
  Let all the assumptions  of Theorem \ref{main} hold.
  For $d$ sufficiently small there exist $C^1$ functions $t\mapsto\bigl(\omega(t),\gamma(t)\bigr)$
  defined for $t\ge 0$ such that the solution $\psi(x,t)$  of  (\ref{SV})
  can be written as in (\ref{sol}-\ref{req})
  with (\ref{omega}-\ref{gamma})  satisfied, and there exists a number
  $\ov M>0$,
  depending only on the initial data, such that
  \be\la{ovY}
  M(T)=\sup\limits_{0\le t\le T}[(1+t)^{3/2}\Vert\chi(t)\Vert_{l^\infty_{-\beta}}
  +(1+t)^3\bigl(|\dot\gamma|+|\dot\omega|\bigr)|]\le \ov M,
  \ee
  uniformly in $T>0$, and $\ov M=O(d)$ as $d\to 0$.
\end{theorem}

\begin{remarks}

(0) This theorem will be deduced from Proposition \ref{ind-arg} in the next section.

  (i)  Theorem \ref{locex} implies that the norms in the definition
  of $M$ are continuous functions of time (and so $M$ is also).

  (ii) The result holds also for negative time with appropriate changes since
  $\psi(x,t)$ solves (\re{S}) if and only if $\overline\psi(x,-t)$ does.

  (iii) The result implies in particular that
  $t^3|\dot\theta-\om|+t^3|\dot\om|\le C$, hence $\om(t)$ and
  $\theta(t)-t\om_{+}$ should converge as $t\to\infty$ while
  $\psi(x,t)-e^{j\theta(t)}\Phi_{\om(t)}(x)$ have limit zero
  in $l^\infty_{-\beta}({\R})$.
\end{remarks}

 %%%%%%%%%%%%%%%%%%%%%%%%%%%%%%%%%%%%%%%%%%%%%%%%%%%%%%%%%%%%%%%%%%%%%%%%%%%%%%%%%%%%%%%%
 %%%%%%%%%%%%%%%%%%%%%%%%%%%%%%%%%%%%%%%%%%%%%%%%%%%%%%%%%%%%%%%%%%%%%%%%%%%%%%%%%%%%%%%%%%
\setcounter{equation}{0}
\section{Proof of transversal decay}
\label{prsec}
%%%%%%%%%%%%%%%%%%%%%%%%%%%%%%%%%%%%%%%%%%%%%%%%%%%%%%%%%%%%%%%%%%%%%%%%%%%%%%%%%%%%%%%%%%
\subsection{Inductive argument (proof of Theorem \ref{yest})}
Let us write the initial data  in the form
  \be\la{0ansatz}
    \psi_0(x)=e^{j\theta_0}(\Phi_{\om_0}(x)+\chi_0(x)).
  \ee
  with $d=\Vert\chi_0\Vert_{l^1_\beta\cap H^1}$ sufficiently small.
By Lemma \ref{in-cond} we can assume that  ${\bf P}^0(\om_0)(\chi_0)=0$
without loss of generality. Then the local existence asserted in Corollary \ref{ME} implies
the existence of an interval $[0,t_1]$ on which are defined $C^1$ functions
$t\mapsto\bigl(\omega(t),\gamma(t)\bigr)$ satisfying (\ref{omega})-(\ref{gamma})
and such that $M(t_1)=\rho$ for some $t_1>0$ and $\rho>0$. By continuity we can make $\rho$
as small as we like by making $d$ and $t_1$ small. The following Proposition
  is proved in section \ref{ptsec} below.
\begin{pro}\la{ind-arg}
In the situation of Theorem \ref{yest}
let $M(t_1)\le \rho$ for some $t_1>0$ and $\rho>0$. Then
there exist numbers
$d_1$ and $\rho_1$, independent of $t_1$, such that
\be\la{half}
M(t_1)\le \rho/2
\ee
if $d=\Vert\chi_0\Vert_{L^1_\beta\cap H^1}<d_1$ and $\rho<\rho_1$.
\end{pro}
Assuming the truth of  Proposition \ref{ind-arg}
for now Theorem
\ref{yest}
will follow from the next argument:\\
Consider the set ${\cal T}$ of $t_1\ge 0$ such that
$\bigl(\omega(t),\gamma(t)\bigr)$ are defined on $[0,t_1]$ and $M(t_1)\le\rho$.
This set is relatively closed by continuity. On the other hand,
(\ref{half}) and Corollary \ref{ME} with sufficiently small $\rho$ and $d$
 imply that
this set is also relatively open, and hence $\sup {\cal
  T}=+\infty$, completing the proof of Theorem \ref{yest}.\bo
%%%%%%%%%%%%%%%%%%%%%%%%%%%%%%%%%%%%%%%%%%%%%%%%%%%%%%%%%%%%%%%%%%%%%%%%%%%%%%%%%%%%%%%%%%%%
\subsection{Frozen linearized equation}
\label{autsec}
%%%%%%%%%%%%%%%%%%%%%%%%%%%%%%%%%%%%%%%%%%%%%%%%%%%%%%%%%%%%%%%%%%%%%%%%%%%%%%%%%%%%%%%%%%%%

A crucial part of the proof of Proposition \ref{ind-arg} is the estimation
of the first term in $M$, for which purpose it is necessary to make
use of the dispersive properties obtained in sections \ref{subspace} and
\ref{rpsec}. Rather
than study directly (\ref{SV2}), whose linear part is non-autonomous,
it is convenient (following \cite{BP,BS}) to introduce  a small
modification of (\ref{sol}), which leads to an autonomous linearized equation.
This new ansatz for the solution is
\be\la{ansatz2}
  \psi(x,t)=e^{j\theta}(\Phi_\om(x)+e^{-j(\theta-\tilde\theta)}\eta),\quad\hbox{where}\;
  \tilde\theta(t)=\om_1t+\theta_0,\;\theta_0=\theta_0\;\hbox{and}\,
  \om_1=\om(t_1)
\ee
so that, $\eta=e^{j(\theta-\tilde\theta)}\chi$ and
$\chi=e^{-j(\theta-\tilde\theta)}\eta$.
Since
$$\dot\chi=e^{-j(\theta-\tilde\theta)}\Bigl(\dot\eta-j(\om+\dot\gamma-\om_1)\eta\Bigr)$$
equation (\ref{SV2}) implies
\be\la{eqeta}
  \dot\eta=j^{-1}(\om_1-\om)\eta+e^{j(\theta-\tilde\theta)}
  {\bf C}\Bigl(e^{-j(\theta-\tilde\theta)}\eta\Bigr)
  +e^{j(\theta-\tilde\theta)}\Bigl(j^{-1}\dot\gamma\Phi_\om-\dot\om\partial_\om\Phi_\om
  +{\bf Q}[e^{-j(\theta-\tilde\theta)}\eta]\Bigr).
\ee
The matrices ${\bf C}$ and $e^{j\phi}$, where $\phi=\theta-\tilde\theta$, do not commute:
\be\la{dif1}
  {\bf C}e^{j\phi}-e^{j\phi}{\bf C}=\delta(x)b\sin\phi~\sigma,\;{\rm where}\;\sigma=
  \left(\ba{cc}
  1    &    0\\
   0   &   -1
   \ea\right),\;b=2a'C^2.
\ee
Using (\ref{dif1}) we rewrite equation (\ref{eqeta}) as
$$
  \dot\eta=j^{-1}(\om_1-\om)\eta+{\bf C}\eta+e^{j(\theta-\tilde\theta)}
  \Bigl(-\delta(x)b\sin(\theta-\tilde\theta)\sigma\eta+j^{-1}\dot\gamma\Phi_\om
  -\dot\om\partial_\om\Phi_\om+{\bf Q}[e^{-j(\theta-\tilde\theta)}\eta]\Bigr).
$$
To obtain a perturbed {\it autonomous} equation we rewrite the first two terms
on the RHS by freezing the coefficients at $t=t_1$. Note that
$$j^{-1}(\omega_1-\omega)+{\bf C}={\bf C_1}-j^{-1}\delta(x)(V-V_1),$$
where $V=a+bP_1$, $V_1=V(t_1)$, and ${\bf C_1}={\bf C}(t_1)$.
The equation for $\eta$ now reads
\be\la{eta-eq}
  \dot\eta={\bf C_1}\eta-j^{-1}\delta(x)(V-V_1)\eta+e^{j(\theta-\tilde\theta)}
  \Bigl(-\delta(x)b\sin(\theta-\tilde\theta)\sigma\eta+j^{-1}\dot\gamma\Phi
  -\dot\om\partial_\om\Phi_\om+{\bf Q}[e^{-j(\theta-\tilde\theta)}\eta]\Bigr)
\ee
The first term is now independent of $t$; the idea is that if
there is sufficiently rapid convergence of $\om(t)$ as $t\to\infty$
the other remaining  terms are small {\it uniformly with respect to} $t_1$.
Finally the equation (\ref{eta-eq}) can be written in the following
{\it frozen form}
\be\la{SV3}
  \dot\eta={\bf C_1}\eta+{\bf f_1}
\ee
where
\be\la{f1}
  {\bf f_1}= -j^{-1}\delta(x)(V-V_1)\eta+e^{j(\theta-\tilde\theta)}
  \Bigl(-\delta(x)b\sin(\theta-\tilde\theta)\sigma\eta+j^{-1}\dot\gamma\Phi
  -\dot\om\partial_\om\Phi_\om+{\bf Q}[e^{-j(\theta-\tilde\theta)}\eta]\Bigr)
\ee
\begin{remark}
  The advantage of (\ref{SV3}) over (\ref{SV2}) is that it can be treated as a perturbed
  autonomous linear equation, so that the estimates from section \ref{subspace} can be used
  directly. The  additional terms in ${\bf f_1}$
  can be estimated as small uniformly in $t_1$: see lemma \ref{alest} below.
  This is the reason for introduction of the  ansatz (\ref{ansatz2}).
\end{remark}
%%%%%%%%%%%%%%%%%%%%%%%%%%%%%%%%%%%%%%%%%%%%%%%%%%%%%%%%%%%%%%%%%%%%%%%%%%%%%%%%%%%%%%%%
\begin{lemma}
\la{alest}
In the situation of Proposition \ref{ind-arg} there exists $c>0$,
independent of $t_1$, such that for $0\le t\le t_1$
$$
|a(t)-a_1|+|b(t)-b_1|+|\theta(t)-\tilde\theta(t)|\le c\rho,
$$
where
\be\la{cor-prop}
\rho:=\sup\limits_{0\le t\le t_1}(1+t^3)(|\dot\gamma(t)|+|\dot\om(t)|)\le M(t_1).
\ee
\end{lemma}
\Pr
By (\ref{cor-prop}), we have
$$
|a(t)-a(t_1)|=|\int\limits_t^{t_1}\dot a(\tau)d\tau|\le c\Bigl(\sup\limits_{0\le\tau\le t_1}
(1+\tau^2)|\dot\om(\tau)|\Bigr)\int\limits_t^{t_1}\fr{d\tau}{1+\tau^2}\le c\rho,
$$
since $|\dot a(\tau)|\le c|\dot\om(\tau)|$.
The difference $|b(t)-b(t_1)|$ can be estimated similarly.
Next
\beqn\nonumber
\theta(t)-\tilde\theta(t)&=&\int_0^t\om(\tau)d\tau+\gamma(t)-\om(t_1)t-\gamma(0)=
\int_0^t(\om(\tau)-\om(t_1))d\tau + \int_0^t\dot\gamma(\tau)d\tau\\
\la{theta-dif}
&=&-\int_0^t\int_{\tau}^{t_1}\dot\om(s)dsd\tau + \int_0^t\dot\gamma(\tau)d\tau.
\eeqn
By (\ref{cor-prop}) the first summand in RHS of (\ref{theta-dif}) can be estimated as
$$
\int_0^t\int_{t_1}^{\tau}|\dot\om(s)|ds\,d\tau\le
\int_0^t\int_{\tau}^{t_1}(1+s)^{2+\ve}|\dot\om(s)|\fr1{(1+s)^{2+\ve}}ds\,d\tau
$$
$$
\le c\sup\limits_{0\le s\le t_1}(1+s)^{2+\ve}|\dot\om(s)|
\int_0^t\int_{\tau}^{t_1}\fr1{(1+s)^{2+\ve}}ds\,d\tau\le c\rho
$$
since the last integral is bounded for $t\in[0,t_1]$.
Finally, for the second summand  on the RHS of (\ref{theta-dif}) inequality (\ref{cor-prop})
implies
$$
|\int_0^t\dot\gamma(\tau)d\tau|\le c\sup\limits_{0\le\tau\le t_1}
(1+\tau^2)|\dot\gamma(\tau)|\int\limits_t^{t_1}\fr{d\tau}{1+\tau^2}\le c\rho
$$
 \bo
%%%%%%%%%%%%%%%%%%%%%%%%%%%%%%%%%%%%%%%%%%%%%%%%%%%%%%%%%%%%%%%%%%%%%%%%%%%%%%%%%%%%%%%%%%

\subsection{Projection onto discrete and continuous spectral spaces}
\label{decompsec}
%%%%%%%%%%%%%%%%%%%%%%%%%%%%%%%%%%%%%%%%%%%%%%%%%%%%%%%%%%%%%%%%%%%%%%%%%%%%%%%%%%%%%%%%%%%%

From sections \ref{subspace} and \ref{rpsec} we have information concerning
$U(t)=e^{{\bf C_1}t}$, in particular decay on the subspace orthogonal
to the (two dimensional) generalized null space. It is
therefore necessary to introduce a further decomposition to take
advantage of this. Recall, by
comparing (\ref{sol}) and (\ref{ansatz2})  that
\be
\la{rtt}
\eta=e^{j(\theta-\tilde\theta)}\chi\quad\hbox{and}\;\;\;{\bf P}^0_t\chi(t)=0.
\ee
Introduce the symplectic projections ${\bf P}^0_1={\bf P}^0_{t_1}$ and
${\bf P}^c_1={\bf P}^c_{t_1}$ onto the discrete and continuous
spectral subspaces defined by the operator ${\bf C_1}$ and write, at
each time $t\in[0,t_1]$:
\be\la{sdec2}
  \eta(t)=g(t)+h(t)
\ee
with $g(t)={\bf P}^0_1\eta(t)$ and  $h(t)={\bf P}^c_1\eta(t)$. The
following lemma shows that it is only necessary to estimate $h(t)$.
\begin{lemma}\label{g}
  In the situation of Proposition \ref{ind-arg}, assume
  $$
   \sup\limits_{0\le t\le t_1}\bigl(
   |\omega(t)-\omega_{1}|+|\theta(t)-\theta_{1}(t)|\bigr)=\Delta
  $$
  is sufficiently small. Then for $0\le t\le t_1$ there exists
  $c(\Delta,\omega_1)$ such that
\be\la{chi}
  c(\Delta,\omega_1)^{-1}\Vert h\Vert _{l^{\infty}_{-\beta}\cap l^2}
  \le\Vert \eta\Vert_{l^{\infty}_{-\beta}\cap l^2}
  \le c(\Delta,\omega_1)\Vert h\Vert _{l^{\infty}_{-\beta}\cap l^2}.
\ee
\end{lemma}
\Pr
This can be proved as in \ci[Lemma 11.5]{BKKS}
\bo
%%%%%%%%%%%%%%%%%%%%%%%%%%%%%%%%%%%%%%%%%%%%%%%%%%%%%%%%%%%%%%%%%%%%%%%%%%%%%%%%%%%%%%%%%%
\subsection{Proof of Proposition \ref{ind-arg}}
\label{ptsec}
%%%%%%%%%%%%%%%%%%%%%%%%%%%%%%%%%%%%%%%%%%%%%%%%%%%%%%%%%%%%%%%%%%%%%%%%%%%%%%%%%%%%%%%%%%%%
To prove Proposition \ref{ind-arg}
we explain how to estimate both terms in $M$,
(\re{ovY}),
 to be $\le \rho/4$,
uniformly in $t_1$.\\
{\it  Estimation of the second term in $M$.}
As in Corollary \ref{ME} we have
$$
  |\dot\gamma(t)|+|\dot\om(t)|\le c_0|\chi(0,t)|^2\le c_0\fr{M(t)^2}{(1+|t|)^3},
  \quad t\le t_1,
$$
since $|\chi(0,t)|\le\Vert\chi(t)\Vert_{l^\infty_{-\beta}}$.
Finally let $\rho_1<1/(4c_0)$ to complete the estimate for the second term in $M$ as
$\le\rho/4$.\\
\noindent
{\it  Estimation of the first term in $M$.}
By Lemma \ref{g} it is enough  to estimate $h$.
Let us apply the projection ${\bf P}^c_1$ to both sides of  (\ref{SV3}).
Then the equation for $h$ reads
\be\la{h-eq}
  \dot h={\bf C_1}h+{\bf P^c_1f_1}
\ee
Now to estimate $h$ we use the Duhamel representation:
\be\la{Duh-rep}
h(t)=U(t) h(0)+\int_0^t U(t-s){\bf P}^c_1{\bf f_1}(s) ds,\quad t\le t_1.
\ee
with $U(t)=e^{{\bf C_1}t}$ the one parameter group just introduced.
Recall that
$
{\bf P}^0_1 h(t)=0$ for $t\in [0,t_1]$.
Therefore
\be\la{UU}
\Vert U(t)h(0)\Vert _{l^{\infty}_{-\beta}}\le c(1+t)^{-3/2}
\Vert h(0)\Vert_{l^1_{\beta}\cap l^2}\le c(1+t)^{-3/2}
\Vert \eta(0)\Vert_{l^1_{\beta}\cap l^2}.
\ee
by  Theorem \ref{TD} and inequalities (\ref{t-small}) and (\ref{chi}).
Let us estimate the integrand on the
right-hand side of (\ref{Duh-rep}).
We use the representation (\re{f1}) for $\bf f_1$ and
apply   Theorem \ref{TD}, Corollary \ref{ME} and Lemma \ref{alest}  to obtain that
\beqn\la{U1}
\Vert U(t-s){\bf P}^c_1{\bf f_1}\Vert _{l^{\infty}_{-\beta}}
&\le& c(1+t-s)^{-3/2}\Vert {\bf P}^c_1({\bf f_1}(t))\Vert_{l^1_{\beta}}
\nonumber\\
&\le& c(1+t-s)^{-3/2} \biggl(|\eta(0,t)|^2+\rho|\eta(0,t)|\biggr)
\nonumber\\
&\le& c(1+t-s)^{-3/2}\biggl(\Vert\eta(t)\Vert^2_{l^\infty_{-\beta}}
  +\rho\Vert\eta(t)\Vert_{l^\infty_{-\beta}}\biggr),\quad t\le t_1.
\eeqn
Now (\ref{chi}), (\ref{Duh-rep}), (\ref{UU}) and (\ref{U1}) imply
$$
\Vert\eta(t)\Vert_{l^\infty_{-\beta}}\le c(1+t)^{-3/2}\Vert\eta(0)\Vert_{l^1_{\beta}\cap l^2}
+c_1\int\limits_0^t \frac {ds}{(1+t-s)^{3/2}}
\biggl(\Vert\eta(s)\Vert^2_{l^\infty_{-\beta}}
  +\rho\Vert\eta(s)\Vert_{l^\infty_{-\beta}}\biggr)
$$
Multiply by $(1+t)^{3/2}$ to deduce
\beqn
(1+t)^{3/2}\Vert\eta(t)\Vert_{l_{-\beta}^\infty}\le cd
&+&c_1\int\limits_{0}^t\frac{(1+t)^{3/2}(1+s)^{-3}}{(1+t-s)^{3/2}}
(1+s)^{3}\Vert\eta(s)\Vert^2_{l^\infty_{-\beta}}ds \\
\nonumber
&+&c_1\rho\int\limits_{0}^t\frac{(1+t)^{3/2}(1+s)^{-3/2}}{(1+t-s)^{3/2}}
(1+s)^{3/2}\Vert\eta(s)\Vert_{l^\infty_{-\beta}}ds
\eeqn
since $\Vert\eta(0)\Vert_{l^1_\beta\cap l^2}\le d$.
Introduce the majorant
$$
m(t):=\sup_{[0,t]}(1+s)^{3/2}\Vert\eta(s)\Vert_{l_{-\beta}^\infty},\quad t\le t_1
$$
and hence
\beqn\label{M}
m(t)\le cd+c_1 m^2(t)\int\limits_{0}^t
\frac{(1+t)^{3/2}(1+s)^{-3}}{(1+t-s)^{3/2}}\,ds
+\rho c_1m(t)\int\limits_{0}^t
\frac{(1+t)^{3/2}(1+s)^{-3/2}}{(1+t-s)^{3/2}}\,ds.
\eeqn
It easy to  see (by splitting up the integrals into
$s<t/2$ and $s\ge t/2$) that both these integrals are bounded independent of $t$.
Thus (\ref{M}) implies that there exist $c,c_2,c_3$, independent of
$t_1$, such that
$$m(t)\le cd+\rho c_2m(t)+c_3m^2(t),\quad t\le t_1.$$
Recall that $m(t_1)\le\rho\le\rho_1$ by assumption. Therefore
this inequality implies that $m(t)$ is bounded for $t\le t_1$, and moreover,
$$ m(t)\le c_4d,\quad t\le t_1$$
if $d$ and $\rho$ are sufficiently small. The constant $c_4$ does
not depend on $t_1$. We choose $d$ in (\ref{close}) small enough that $d<\rho/(4c_4)$.
Therefore,
$$\sup_{[0,t_1]}(1+t)^{3/2}\Vert\eta(t)\Vert_{l_{-\beta}^\infty} < \rho/4$$
if $d$ and $\rho$ are sufficiently small. This bounds the first
term as $<\rho/4$ by (\ref{rtt}) and hence
$M(t_1)< \rho/2$, completing the proof of Proposition \ref{ind-arg}.
\bo
%%%%%%%%%%%%%%%%%%%%%%%%%%%%%%%%%%%%%%%%%%%%%%%%%%%%%%%%%%%%%%%%%%%%%%%%%%%%%%%%%%%%%%%%%%
\setcounter{equation}{0}
\section{Soliton asymptotics}
\label{solas-sec}
%%%%%%%%%%%%%%%%%%%%%%%%%%%%%%%%%%%%%%%%%%%%%%%%%%%%%%%%%%%%%%%%%%%%%%%%%%%%%%%%%%%%%%%%%%%%

Here we prove our main Theorem \ref{main} using the bounds (\ref{ovY}).
For the solution $\psi(x,t)$ to (\ref{S}) let us define the accompanying soliton as
$s(x,t)=\psi_{\om(t)}(x)e^{i\theta(t)}$, where  $\dot\theta(t)=\om(t)+\dot\gamma(t)$.
Then for the difference $z(x,t)=\psi(x,t)-s(x,t)$ we obtain easily from
 equations (\ref{S}) and (\ref{NEP})
\be\la{z}
i\dot z(x,t)=-z''(x,t)+\dot\gamma s(x,t)-i\dot\om\partial_{\om}s(x,t)
-\delta(x)\Bigl(F(\psi(0,t))-F(s(0,t))\Bigr).
\ee
Then
\be\la{z1}
z(t)=W(t)z(0)+\int\limits_0^t W(t-\tau)f(\cdot,\tau)d\tau,
\ee
where $f(x,t)=\dot\gamma s(x,t)-i\dot\om\partial_{\om}s(x,t)
-\delta(x)\Bigl(F(\psi(0,t))-F(s(0,t))\Bigr)$, and $W(t)$ is the dynamical group
of the free Schr\"odinger equation.
Since $\gamma(t)-\gamma_{+}$,
$\om(t)-\om_{+}=\cO(t^{-2})$, and  therefore
$\theta(t)-\om_{+}t-\gamma_{+}=\cO(t^{-1})$ for $t\to\infty$,
to establish the asymptotic behavior (\ref{sol-as}) it suffices to
prove that
\be\la{td}
z(t)=W(t)\Phi_{+}+r_{+}(t)
\ee
with some $\Phi_{+}\in l^2(\Z)$ and $\Vert r_{+}(t)\Vert=\cO(t^{-1/2})$.
Let us rewrite (\ref{z1}) as
\be\la{z2}
z(t)=W(t)\Bigl(z(0)+\int\limits_0^{\infty}W(-\tau)f(\cdot,\tau)d\tau\Bigl)
-\int\limits_0^tW(t-\tau)f(\cdot,\tau)d\tau=W(t)\Phi_{+}+r_{+}(t)  .
\ee
Let us recall that by (\ref{ovY})
$$
|\dot\om(t)|\le c {(1+t)^{-3}},\;|\dot\gamma(t)|\le c {(1+t)^{-3}},\;
|F(\psi(0,t))-F(s(0,t))|\le c|\chi(0,t)|\le c {(1+t)^{-3/2}}.
$$
Hence, the unitarity in $l^2$ of the group $W(t)$ implies that
 $\Phi_{+}=z(0)+\ds\int\limits_0^{\infty}W(-\tau)f(\cdot,\tau)d\tau\in l^2$, and
 $\Vert r_+(t)\Vert=\cO(t^{-1/2}),\;t\to\infty$.
\bo
\setcounter{section}{0}
\setcounter{equation}{0}
\protect\renewcommand{\thesection}{\Alph{section}}
\protect\renewcommand{\theequation}{\thesection. \arabic{equation}}
\protect\renewcommand{\thesubsection}{\thesection. \arabic{subsection}}
\protect\renewcommand{\thetheorem}{\Alph{section}.\arabic{theorem}}

%%%%%%%%%%%%%%%%%%%%%%%%%%%%%%%%%%%%%%%%%%%%%%%%%%%%%%%%%%%%%%%%%%%%%%%%%%%%%%%%%%%%%%%%%
 %%%%%%%%%%%%%%%%%%%%%%%%%%%%%%%%%%%%%%%%%%%%%%%%%%%%%%%%%%%%%%%%%%%%%%%%%%
\section{The resolvent}
%\la{appendix}
%%%%%%%%%%%%%%%%%%%%%%%%%%%%%%%%%%%%%%%%%%%%%%%%%%%%%%%%%%%%%%%%%%%%%%%%%%%%%%
\subsection*{A.1 Calculation of the matrix kernel}
\label{res}
%%%%%%%%%%%%%%%%%%%%%%%%%%%%%%%%%%%%%%%%%%%%%%%%%%%%%%%%%%%%%%%%%%%%%%%%%%%%%%%%
Here we will construct  matrix  kernel of the resolvent ${\bf R}(\lam)$ explicitly
\be\la{mic}
   {\bf R}(\lam,x,y)= \left( \ba{cc}R_{11}(\lam,x,y)&R_{12}(\lam,x,y)\\
   R_{21}(\lam,x,y)&R_{22}(\lam,x,y) \ea\right)
\ee
which is the solution to the equation
\be\la{mice}
   ({\bf C}-\lam){\bf R}(\lam,x,y)=\de(x-y)
   \left( \ba{cc}
   1&0\\
   0&1
   \ea\right).
   \ee
%%%%%%%%%%%%%%%%%%%%%%%%%%%%%%%%%%%%%%%%%%%%%%%%%%%%%%%%%%%%%%%%%%%%%%%%%%%%%%%%%%%%%%%%%%
%%%%%%%%%%%%%%%%%%%%%%%%%%%%%%%%%%%%%%%%%%%%%%%%%%%%%%%%%%%%%%%%%%%%%%%%%%%%%%%%%%%%%%%%%
{\bf Calculation of first column}
For the first column
$R_I(\lam,x,y):= \left(\ba{c}R_{11}(\lam,x,y)\\
R_{21}(\lam,x,y) \ea\right)$  of the matrix ${\bf R}(\lam,x,y)$ we obtain
\be\la{mice1}
   ({\bf C}-\lam)R_I(\lam,x,y)=\de(x-y)
    \left( \ba{c}
   1\\
   0 \ea\right).
   \ee
If $x\ne 0$ and
$x\ne y$,  (\ref{mice1}) takes the form
\be\la{hom}
   \left(
   \ba{lcl}-\lam    &&        {\bf D}_2\\
   -{\bf D}_1       &&        -\lam \ea \right)R_I(\lam,x,y) =\left(
   \ba{rcr}-\lam    &&        -\Delta_L-\om\\
   \Delta_L+\om     &&        -\lam \ea \right)R_I(\lam,x,y) =
   0,~~~~~x\ne 0,~~x\ne y.
\ee
The general solution is a linear combination of exponential solutions of type $e^{i\theta x}v$.
Substituting into (\re{hom}), we get
\be\la{homh}
   \left(
   \ba{rcr}-\lam            &&    2-2\cos\theta+\om\\
   -2+2\cos\theta-\om       &&    -\lam
   \ea
   \right)v=0.
\ee
For nonzero vectors $v$, the determinant of the matrix vanishes,
$$
   \lam^2+(2-2\cos\theta+\om)^2=0.
$$
Finally, we obtain four roots $\pm\theta_\pm(\lam)$ in
$D:=-\pi\le\rRe\theta\le\pi$ with
\be\la{ik}
 2- 2\cos\theta_\pm(\lam)=-\om\mp i\lam.
\ee
We choose the cuts in the complex plane $\lam$:\\
the cut ${\cal C}_+:=[i\om,i(\om+4)]$ for $\theta_+(\lam)$
and the cut ${\cal C}_-:=[-i(\om+4),-i\om]$ for $\theta_-(\lam)$ and
\be\la{re}
   \rIm\theta_\pm(\lam)>0,~~~~~~~\lam\in\C\setminus {\cal C}_\pm.
\ee
It remains to derive the vector $v=(v_1,v_2)$ which is solution to (\re{homh}):
$$
   v_2=-\fr{2-2\cos\theta_\pm+\om}{\lam}v_1=\fr{ \pm i\lam}{\lam}v_1=\pm i v_1.
$$
Therefore, we have two corresponding vectors
$v_\pm=\left(
\ba{r}1\\
\pm i
\ea
\right)$
and we get four linearly independent exponential solutions
$v_+e^{\pm i\theta_+ x}$ and  $v_-e^{\pm i\theta_- x}$.

Now we can solve the equation (\re{mice1}).
First we rewrite it using the representation (\re{lin4})
for the operator ${\bf C}$,
\be\la{homr}
   \!\left(\!
   \ba{rr}-\lam          &     -\Delta_L+\om\\
   \Delta_L-\om          &     -\lam
   \ea
   \!\right)\!\left(\ba{c}
   \!R_{11}(\lam,x,y)\\
   R_{21}(\lam,x,y) \ea\right) =\de(x-y) \left( \ba{c}
   1\\
   0
   \ea\!\right)+
   \de(x)
   \left(
   \!\ba{cc}
   0      &a\\
   -a-b &0
   \ea\!\right)\!\left(\ba{c}
   R_{11}(\lam,0,y)\\
   R_{21}(\lam,0,y)
   \ea\right)
\ee
%%%%%%%%%%%%%%%%%%%%%%%%%%%%%%%%%%%%%%%%%%%%%%%%%%%%%%%%%%%%%%%%%%%%%%%%%%%%%%%%%%%%%%%%%5
Let us consider $y>0$ for the concreteness. Then
the RHS vanishes in the open intervals $(-\infty,0)$,$(0,y)$ and $(y,\infty)$.
Hence, for the parameter $\lam$ outside the cuts $\cC_\pm$,
the solution admits the representation
\be\la{ABC}
   R_{I}(\lam,x,y)
   =\left\{ \ba{ll}
   A_{+}e^{-i\theta_{+}x}v_{+} +A_{-}e^{-i\theta_{-} x}v_{-},  \!&\!x<0,\\\\
   B^{-}_{+}e^{-i\theta_{+} x}v_{+}+B^{-}_{-}e^{-i\theta_{-} x}v_{-}
   +B^{+}_{+}e^{i\theta_{+} x}v_{+}+B^{+}_{-}e^{i\theta_{-} x}v_{-},  \!&\!0<x<y,\\\\
   C_{+}e^{i\theta_{+} x}v_{+}+C_{-}e^{i\theta_- x}v_{-},  \!&\!x>y
   \ea
   \right.
\ee
since by (\re{re}), the exponent $e^{-i\theta_\pm x}$ decays for $x\to -\infty$,
and similarly, $e^{i\theta_\pm x}$ decays for $x\to\infty$.
Next we need eight equations to calculate the eight constants $A_+,\dots,C_-$.
We have two continuity equations and two jump conditions for the derivatives
at the points $x=0$ and $x=y$. These four vector equations
give just eight scalar equations for the calculation.
%%%%%%%%%%%%%%%%%%%%%%%%%%%%%%%%%%%%%%%%%%%%%%%%%%%%%%%%%%%%%%%%%%%%%%%%%%%%%%%
\\
{\bf Continuity at $x=y$:} $R_I(y-0,y)=R_I(y+0,y)$, i.e.
$$
  B^{-}_{-}v_{+}/e_{+}+B^{-}_{-}v_{-}/e_{-}+B^{+}_{+}v_{+}e_{+}+B^{+}_{-}v_{-}e_{-}
  =C_{+}v_{+}e_{+}+C_{-}v_{-}e_{-},
$$
where  $e_\pm:= e^{i\theta_\pm y}$.
It is equivalent to
\be\la{xy}
   \left\{
   \ba{l}
   B^{-}_{+}/e_{+}+B^{+}_{+}e_{+}=C_{+}e_{+},\\\\
   B^{-}_{-}/e_{-}+B^{+}_{-}e_{-}=C_{-}e_{-}.
   \ea
   \right.
\ee
{\bf Continuity at $x=0$:} $R_I(-0,y)=R_I(+0,y)$, i.e.
$$
   A_{+}v_{+}+A_{-}v_{-}=B^{-}_{+}v_{+}+B^{-}_{-}v_{-}+B^{+}_{+}v_{+}+B^{+}_{-}v_{-}
$$
that is equivalent to
\be\la{x0}
   \left\{
   \ba{l}
   A_{+}=B^{-}_{+}+B^{+}_{+},\\\\
   A_{-}=B^{-}_{-}+B^{+}_{-}.
   \ea
   \right.
\ee
%%%%%%%%%%%%%%%%%%%%%%%%%%%%%%%%%%%%%%%%%%%%%%%%%%%%%%%%%%%%%%%%%%%%%%%%%%%%%%%%%%%%%%%%%%%
{\bf ``Jump" at $x=y$:} Noting that
\be\la{vpm}
   \left(
   \ba{r}
   1\\
   0
   \ea\right)=\ds\fr{v_{+}+v_{-}}2.
\ee
At $x=y$ equation (\ref{homr}) reads
\be\la{x0y}
   \left\{
   \ba{l}
   -\lam C_{+}e_{+}-iB_{+}^{-}e^{i\theta_{+}}/e_{+}-iB_{+}^{+}e^{-i\theta_{+}}e_{+}
   -iC_{+}e_{+}e^{i\theta_{+}}+(2+\om)iC_{+}e_{+}=1/2,\\\\
   -\lam C_{-}e_{-}+iB_{-}^{-}e^{i\theta_{-}}/e_{-}+iB_{-}^{+}e^{-i\theta_{}}e_{-}
   +iC_{-}e_{-}e^{i\theta_{-}}-(2+\om)iC_{-}e_{-}=1/2.
   \ea
   \right.
\ee
Also we need equation (\ref{homr}) at $x=y+1$:
$$
   \left\{
   \ba{l}
   -\lam C_{+}e^{i\theta_{+}(y+1)}-iC_{+}e^{i\theta_{+}y}-iC_{+}e^{i\theta_{+}(y+2)}
   +(2+\om)iC_{+}e^{i\theta_{+}(y+1)}=0,\\\\
   -\lam C_{-}e^{i\theta_{-}(y+1)}+iC_{-}e^{i\theta_{-}y}+iC_{-}e^{i\theta_{-}(y+2)}
   -(2+\om)iC_{-}e^{i\theta_{-}(y+1)}=0.
   \ea
   \right.
$$
Hence
\be\la{sub}
   \left\{
   \ba{l}
   -\lam C_{+}-iC_{+}e^{i\theta_{+}}
   +(2+\om)iC_{+}=iC_{+}e^{-i\theta_{+}},\\\\
   -\lam C_{-}+iC_{-}e^{i\theta_{-}}
   -(2+\om)iC_{-}=-iC_{-}e^{-i\theta_{-}}.
   \ea
   \right.
\ee
Combining  (\ref{x0y}) and (\ref{sub}), we get
$$
   \left\{
   \ba{l}
   -iB_{+}^{-}e^{i\theta_{+}}/e_{+}-iB_{+}^{+}e^{-i\theta_{+}}e_{+}
   +iC_{+}e_{+}e^{-i\theta_{+}}=1/2,\\\\
   iB_{-}^{-}e^{i\theta_{-}}/e_{-}+iB_{-}^{+}e^{-i\theta_{}}e_{-}
   -iC_{-}e_{-}e^{-i\theta_{-}}=1/2.
   \ea
   \right.
$$
After substituting of $C_{\pm}$ from (\re{xy}), the constants
$B^{+}_{\pm}$ cancel and we get
$$
B_{+}^{-}\frac{i}{e_{+}}(-e^{i\theta_{+}}+e^{-i\theta_{+}})=1/2,\quad
B_{-}^{-}\frac{i}{e_{-}}(e^{i\theta_{-}}-e^{-i\theta_{-}})=1/2
$$
and then
\be\la{BB}
B^{-}_{+}=\ds\fr {e_+}{4\sin\theta_+},~~~~~~~~~~~~~~~~B^{-}_{-}=-\ds\fr{e_-}{4\sin\theta_-}.
\ee
%%%%%%%%%%%%%%%%%%%%%%%%%%%%%%%%%%%%%%%%%%%%%%%%%%%%%%%%%%%%%%%%%%%%%%%%%%%%%%%%%%%%%%%%55
{\bf ``Jump" at $x=0$:} Noting that
$$
\left(\ba{cc}0&a\\-a-b&0\ea\right)v_{+}=\left(\ba{c}ia\\-a-b\ea\right)
=iv_{+}(a+\frac b2)-iv_{-}\frac b2=iv_{+}\al-iv_{-}\beta,
$$
$$
\left(\ba{cc}0&a\\-a-b&0\ea\right)v_{-}=\left(\ba{c}-ia\\-a-b\ea\right)
=iv_{+}\frac b2-iv_{-}(a+\frac b2)=iv_{+}\beta-iv_{-}\al,
$$
where
$\al=a+\ds\fr{b}2$, $\beta=\ds\fr{b}2$.
Hence,
$$
\left(\ba{cc}0&a\\-a-b&0\ea\right)\!(A_{+}v_{+}+A_{-}v_{-})
=iv_{+}(A_{+}\al+A_{-}\beta)-iv_{-}(A_{+}\beta+A_{-}\al)
$$
At $x=0$ equation (\ref{homr}) reads
\be\la{x00}
   \left\{
   \ba{l}
   -\lam A_{+}-iA_{+}e^{i\theta_{+}}-iB_{+}^{-}e^{-i\theta_{+}}
   -iB_{+}^{+}e^{i\theta_{+}}+(2+\om)iA_{+}=i(A_{+}\al+A_{-}\beta),\\\\
   -\lam A_{-}+iA_{-}e^{i\theta_{-}}+iB_{-}^{-}e^{-i\theta_{+}}
   +iB_{-}^{+}e^{i\theta_{-}}-(2+\om)iA_{-}=-i(A_{+}\beta+A_{-}\al)
   \ea
   \right.
\ee
Also we need equation (\ref{homr}) at $x=y+1$:
$$
\left\{
   \ba{l}
   -\lam A_{+}e^{i\theta_{+}}-iA_{+}-iA_{+}e^{2i\theta_{+}}+(2+\om)iA_{+}e^{i\theta_{+}}=0,\\\\
   -\lam A_{-}e^{i\theta_{-}}+iA_{-}+iA_{-}e^{2i\theta_{-}}-(2+\om)iA_{-}e^{i\theta_{-}}=0.
   \ea
   \right.
$$
Therefore
\be\la{x01}
\left\{
   \ba{l}
   -\lam A_{+}-iA_{+}e^{-i\theta_{+}}-iA_{+}e^{i\theta_{+}}+(2+\om)iA_{+}=0,\\\\
   -\lam A_{-}+iA_{-}e^{-i\theta_{-}}+iA_{-}e^{i\theta_{-}}-(2+\om)iA_{-}=0.
   \ea
   \right.
\ee
Substituting (\ref{x01}) into (\ref{x00}) we get
$$
\left\{
   \ba{l}
   -iB_{+}^{-}e^{-i\theta_{+}}-iB_{+}^{+}e^{i\theta_{+}}+iA_{+}e^{-i\theta_+}
   =i(A_{+}\al+A_{-}\beta),\\\\
   iB_{-}^{-}e^{-i\theta_{+}}+iB_{-}^{+}e^{i\theta_{-}}-iA_{-}e^{-i\theta_-}
   =-i(A_{+}\beta+A_{-}\al)
   \ea
   \right.
$$
Substituting here (\re{x0}), we get after cancellations,
$$
\left\{
\ba{l}
(e^{i\theta_+}-e^{-i\theta_+}+\al)B^{+}_{+}+\beta B^{+}_{-}=-\al B^{-}_{+}-\beta B^{-}_{-}
\\
\\
\beta B^{+}_{+}+(e^{i\theta_-}-e^{-i\theta_-}+\al)B^{+}_{-}=-\beta B^{-}_{+}-\al B^{-}_{-}
\ea\right.
$$
Hence, the solution is given by
\be\la{BBB}
\left(
\ba{r}
B^{+}_{+}\\
B^{+}_{-}
\ea\right)
=-\fr 1{D}
\left(
\ba{cc}
2i\sin\theta_{-}+\al&  -\beta\\
-\beta     &  2i\sin\theta_{+}+\al
\ea\right)
\left(
\ba{cc}
\al&\beta\\
\beta&\al
\ea\right)
\left(
\ba{r}
B^{-}_{+}\\
B^{-}_{-}
\ea\right),
\ee
where $D$ is the determinant
\be\la{dete}
D:=(2i\sin\theta_++\al)(2i\sin\theta_-+\al)-\beta^2,
\ee
and $B^{-}_{+},B^{-}_{-}$ are given by (\re{BB}).
The formulas (\re{BB}) and (\re{BBB}) imply
\be\la{BBBB}
B^{+}_{+}=\frac 1{2D}\left(-\ds\frac{2i\al\sin\theta_{-}+\al^2-\beta^2}{2\sin\theta_+}e_{+}
+i\beta e_{-}\!\right)\!,\quad
B^{+}_{-}=\frac 1{2D}\left(-i\beta e_{+}
+\ds\frac{2i\al\sin\theta_++\al^2-\beta^2}{2\sin\theta_-}e_{-}\!
 \!\right)
\ee
Using the identities
\beqn\nonumber
2i\al\sin\theta_-+\al^2-\beta^2&=&D-2i\al\sin\theta_++4\sin\theta_+\sin\theta_-,\\
\nonumber
2i\al\sin\theta_++\al^2-\beta^2&=&D-2i\al\sin\theta_-+4\sin\theta_+\sin\theta_-,
\eeqn
let us rewrite (\ref{BBBB}) as
\be\la{fin}
B^{+}_{+}=-\ds\frac{e_+}{4\sin\theta_+}+\frac 1{2D}
\Bigl((i\al-2\sin\theta_{-})e_{+}+i\beta e_{-}\Bigr),
\quad
B^{+}_{-}=\ds\frac{e_-}{4\sin\theta_-}-\frac 1{2D}
\Bigl(i\beta e_{+}+(i\al-2\sin\theta_{+})e_{-}\Bigr).
\ee
Finally, the formulas  (\re{ABC})--(\re{x0}),
(\re{BB}) and (\re{fin}) give the first column $R_I(\lam,x,y)$ of
the resolvent for $y>0$:
\be\la{RI}
R_{I}(\lam,x,y)=\Gamma_{I}(\lam,x,y)+P_{I}(\lam,x,y),
\ee
where
\be\la{GI}
\Gamma_{I}(\lam,x,y)=\frac 1{4\sin\theta_+}(e^{i\theta_{+}|x-y|}
-e^{i\theta_{+}(|x|+|y|)})v_{+}
-\frac 1{4\sin\theta_-}(e^{i\theta_{-}|x-y|}-e^{i\theta_{-}(|x|+|y|)})v_{-},
\ee
and
\beqn\la{PI}
P_{I}(\lam,x,y)&=&\frac 1{2D}\Bigl[\Bigl((i\al-2\sin\theta_{-})e^{i\theta_{+}(|x|+|y|)}
+i\beta e^{i(\theta_{+}|x|+\theta_{-}|y|)}\Bigr)v_{+}\\
\nonumber
&-&\Bigl(i\beta e^{i(\theta_{-}|x|+\theta_{+}|y|)}
+(i\al-2\sin\theta_{+})e^{i\theta_{-}(|x|+|y|)}
\Bigr)v_{-}\Bigr]
\eeqn
\\
%%%%%%%%%%%%%%%%%%%%%%%%%%%%%%%%%%%%%%%%%%%%%%%%%%%%%%%%%%%%%%%%%%%%%%%%%%%%%%%%%%%%%%%%%%%%%
{\bf Calculation of second column}
The second column is given by similar  formulas with the vector
$\left(
\ba{r}
0\\
1
\ea\right)
$
instead of
$\left(
\ba{r}
1\\
0
\ea\right)
$
in (\re{homr}).
Respectively, (\re{vpm}) is changed by
$$
\left(
\ba{r}
0\\
1 \ea\right)=i\ds\fr{v_--v_+}{2}.
$$
Hence, we have now change $1/2$ by $-i/2$ in the first equation of (\re{xy}) and $1/2$
by $i/2$ in the second one. Respectively, (\re{BB}) for the second column reads
$$
B^{-}_{+}=-\ds\fr {ie_+}{4\sin\theta_+},~~~~~~~~~~~~~~~~
B^{-}_{-}=-\ds\fr {ie_-}{4\sin\theta_-}.
$$
Then the second column
$R_{II}(\lam,x,y)$ of the resolvent reads:
\be\la{RII}
R_{II}(\lam,x,y)=\Gamma_{II}(\lam,x,y)+P_{II}(\lam,x,y),
\ee
where
\be\la{GII}
\Gamma_{II}(\lam,x,y)=-\frac i{4\sin\theta_+}(e^{i\theta_{+}|x-y|}
-e^{i\theta_{+}(|x|+|y|)})v_{+}
-\frac i{4\sin\theta_-}(e^{i\theta_{-}|x-y|}-e^{i\theta_{-}(|x|+|y|)})v_{-},
\ee
and
\beqn\la{PII}
P_{II}(\lam,x,y)&=&\frac i{2D}\Bigl[\Bigl(-(i\al-2\sin\theta_{-})e^{i\theta_{+}(|x|+|y|)}
+i\beta e^{i(\theta_{+}|x|+\theta_{-}|y|)}\Bigr)v_{+}\\
\nonumber
&+&\Bigl(i\beta e^{i(\theta_{-}|x|+\theta_{+}|y|)}
-(i\al-2\sin\theta_{+})e^{i\theta_{-}(|x|+|y|)}\Bigr)v_{-}\Bigr]
\eeqn
Note, that if $y<0$ we get the same formulas.
%%%%%%%%%%%%%%%%%%%%%%%%%%%%%%%%%%%%%%%%%%%%%%%%%%%%%%%%%%%%%%%%%%%%%%%%%%%%%%%%%%%%%%%%%%%
%%%%%%%%%%%%%%%%%%%%%%%%%%%%%%%%%%%%%%%%%%%%%%%%%%%%%%%%%%%%%%%%%%%%%%%%%%%%%%%%%%%%%%%%%%%%
\subsection*{A.2 The poles of the resolvent}
\label{polres}
%%%%%%%%%%%%%%%%%%%%%%%%%%%%%%%%%%%%%%%%%%%%%%%%%%%%%%%%%%%%%%%%%%%%%%%%%%%%%%%%%%%%%%%%%
The poles of the resolvent correspond to the roots of the determinant (\re{dete}),
\be\la{deter}
  D(\lam):=\al^2+2i\al(\sin\theta_++\sin\theta_-)
-4\sin\theta_+\sin\theta_--\beta^2=0.
\ee
with $\theta_\pm$ as in (\ref{ik})-(\ref{re}). Thus $D(\lam)$ is an analytic
function on $\C\setminus {\cal C}_-\cup{\cal C}_+$.
Since there are two possible
values for the square roots in $\theta_\pm$ there is a corresponding
four-sheeted function $\tilde D(\lam)$ analytic on a four sheeted
cover
of $\C$ which is branched over ${\cal C}_-$ and ${\cal C}_+$. We call
the sheet defined by (\ref{re}) the {\it physical sheet}.

We will reduce the equation (\re{deter}) to the solution of two
successive quadratic equations. These can be solved explicitly but
the process involves squaring and thus actually produces zeros
of the function $\tilde D(\lam)$ rather than of $D(\lam)$.
Therefore we will then have to check whether or not the roots do
actually lie on the physical sheet.
%%%%%%%%%%%%%%%%%%%%%%%%%%%%%%%%%%%%%%%%%%%%%%%%%%%%%%%%%%%%%%%%%%%%%%%%%%%%%
\begin{pro}\la{lam0}
If $a'=(4a+a^3)/4C^2$ then $\lam=0$ is a root of the determinant $D(\lam)$ with
multiplicity 4, otherwise $\lam=0$ is a root of the determinant
$D(\lam)$ with multiplicity 2.
\end{pro}
%%%%%%%%%%%%%%%%%%%%%%%%%%%%%%%%%%%%%%%%%%%%%%%%%%%%%%%%%%%%%%%%%%%%%%%%%%%%%        \\
\Pr
First let us check that $\lam=0$ is a root of $D(\lam)$. Let us
represent $\om$ by means of $a=a(C^2)$. By results of \S\ref{swsec}
$$
\cosh k=\ds\fr{|\om+2|}2,\quad \sinh k=\ds\fr{|a|}2.
$$
Therefore
\be\la{aom}
(\om+2)^2=4+a^2.
\ee
$$
  D(0)=\alpha^2-\beta^2+2i\alpha ia+a^2=(a+b/2)^2-b^2/4-2(a+b/2)a+a^2=0
$$
since
$$
\sin\theta_{\pm}=\sqrt{1-(\om+2)^2/4}=\sqrt{1-(4+a^2)/4}=\sqrt{-a^2/4}=ia/2.
$$
Now let us compute $D'(\lam)$:
\beqn\nonumber
D'(\lam)&=&i\alpha\Big(\fr{-i(\om+2+i\lam)}{\sqrt{4-(\om+2+i\lam)^2}}
+\fr{i(\om+2-i\lam)}{\sqrt{4-(\om+2-i\lam)^2}}\Big)\\
\nonumber
&+&\fr{i(\om+2+i\lam)\sqrt{4-(\om+2-i\lam)^2}}{\sqrt{4-(\om+2+i\lam)^2}}
-\fr{i(\om+2-i\lam)\sqrt{4-(\om+2+i\lam)^2}}{\sqrt{4-(\om+2-i\lam)^2}}.
\eeqn
Hence $D'(0)=0$ and $\lam=0$ is the root of $D(\lam)$ of
multiplicity at least 2. Further calculation shows that the Taylor
series for $D$ near zero takes the form:
\be
\label{tayld}
D(\lam)=\frac{8a+2a^3-4b}{a^3}\lam^2+O(\lam^4).
\ee
Therefore $\lam=0$ is the root of $D(\lam)$ of multiplicity 4 if and only if
$b=(4a+a^3)/2$, i.e. $a'=(4a+a^3)/4C^2$.
\bo

%%%%%%%%%%%%%%%%%%%%%%%%%%%%%%%%%%%%%%%%%%%%%%%%%%%%%%%%%%%%%%%%%%%%%%%%%%%%%%
Now we prove that there exist a values of $a'C^2=f(a)$ such that
the only root of $D(\lam)$ is $\lam =0$ with multiplicity $2$.
\begin{pro}
  If $a'=-a/C^2$, then the only root of $D(\lam)$
  is $\lam =0$ with multiplicity $2$.
\end{pro}
\Pr
One has $b=-2a$, $\al=a+\fr b2=0$ and then equation (\ref{deter}) implies
\be\la{ss}
\sin\theta_+\sin\theta_-=-\beta^2/4=-a^2/4<0.
\ee
From the other hand, by (\ref{aom}) we get
\beqn\la{sisi}
\sin^2\theta_+\sin^2\theta_-&=&\Big(1-\fr{(\om+2+i\lam)^2}4\Big)
\Big(1-\fr{(\om+2-i\lam)^2}4\Big)\\
\nonumber
&=&\Big(1-\fr{(\om+2)^2}4+\fr{\lam^2}4\Big)^2+\fr{(\om+2)^2\lam^2}4
=\Big(-\fr{a^2}4+\fr{\lam^2}4\Big)^2+\lam^2+\frac{a^2\lam^2}4
\eeqn
Hence, (\ref{ss}) and (\ref{sisi}) imply
$$
\fr{a^4}{16}=\Big(-\fr{a^2}4+\fr{\lam^2}4\Big)^2+\lam^2+\frac{a^2\lam^2}4
$$
Hence
$$
\lam^2(\lam^2+2a^2+16)=0
$$
It remains to check that the roots $\lam=\pm i\sqrt{2a^2+16}$
don't lie on the physical branch.\\
%Using (\ref{ss}) we obtain
%\beqn\nonumber
%(\sin\theta_++\sin\theta_-\!)^2\!\!\!&=&\!\!\!\sin^2\theta_+\!+\sin^2\theta_-\!
%+2\sin\theta_+\sin\theta_-=2-\fr{(\om+2+i\lam)^2}4-\fr{(\om+2-i\lam)^2}4-\fr{a^2}2\\
%\la{sisi1}
%\!\!\!&=&\!\!\!2-\fr{(\om+2)^2}2+\fr{\lam^2}2-\fr{a^2}2
%=-a^2+\fr{\lam^2}2=-2a^2-8<0
%\eeqn
Let $\theta_+=x_++iy_+$, $\theta_-=x_-+iy_-$ with $y_{\pm}>0$.
Then
$$
\cos\theta_+=\cos x_+\cosh y_+-i\sin x_+\sinh y_+=\fr{\om+2+\sqrt{2a^2+16}}2
=\fr{\sqrt{a^2+4}+\sqrt{2a^2+16}}2>0,
$$
$$
\cos\theta_-=\cos x_-\cosh y_--i\sin x_-\sinh y_-=\fr{\om+2-\sqrt{2a^2+16}}2
=\fr{\sqrt{a^2+4}-\sqrt{2a^2+16}}2<0.
$$
Hence, $x_+=0$, $x_-=\pi$ and
$$
\sin\theta_+=\sin(iy_+)=i\sinh y_+,
\quad\sin\theta_-=\sin(\pi+iy_-)=-i\sinh y_-
$$
Finally we obtain that $\sin\theta_+\sin\theta_-=\sinh y_+\sinh y_->0$,
which contradict (\ref{ss}).
\bo

%%%%%%%%%%%%%%%%%%%%%%%%%%%%%%%%%%%%%%%%%%%%%%%%%%%%%%%%%%%%%%%%%%%%%
%%%%%%%%%%%%%%%%%%%%%%%%%%%
%%%%%%%%%%%%%%%%%%%%%%%%%%%%%%%%%%%%%%%%%%%%%%%%%%%%%%%%%%%%%%%%%%%%%
%%%%%%%%%%%%%%%%%%%%%%%%%%%

\end{document}